\documentclass{article}
\usepackage{graphicx}
\usepackage{amsmath}
\usepackage{amssymb}
\usepackage{url}

\usepackage[T1]{fontenc}
\usepackage{mathpazo}
\usepackage{diagrams}
\diagramstyle{midshaft,PostScript=dvips,nohug}
\newarrow{Dotsto}.. ..>

\newtheorem{lemma}{Lemma}[section]
\newtheorem{prop}[lemma]{Proposition}

\newtheorem{cor}[lemma]{Corollary}

\newcommand{\id}{\operatorname{id}}

\newcommand{\Spec}{\operatorname{Spec}}
\newcommand{\A}{\mathcal{A}}
\newcommand{\E}{\mathcal{E}}
\newcommand{\M}{\mathfrak{M}}
\newcommand{\Ss}{\mathcal{S}}
\newcommand{\GG}{\mathcal{G}}

\title{New methods for old spaces:\\ synthetic differential geometry\\
}
\author{Anders Kock}
\date{}

\begin{document}
\maketitle

\medskip

\section*{Introduction}
The synthetic method consists in  consideration of a class of objects 
in terms of their (often axiomatically assumed)  mutual relationship, say 
incidence relations, 
- disregarding what the objects are ``made up of''. For geometry, this 
method  goes back to the time of Euclid. With the advent of category 
theory, it became possible to make  the notion of ``relationship'' more 
precise, in terms of the {\em maps} in  some category $\E$ of 
``spaces'', in some broad sense of this word.  

A synthetic theory  might be presented by  giving axioms 
for some  good category $\E$, possibly with some added structure. Thus, the 
list of axioms could begin: ``let $\E$ be a good category, and let $R$ be a 
commutative ring object (to be thought of as the number line) in $\E$ 
\ldots''. Even though this is the way most texts in synthetic 
differential geometry begin, some texts are purely 
synthetic/combinatorial, presupposing a category $\E$, but do not 
presuppose any  ring object $R$; this applies e.g.\ to  Sections 
\ref{Touchx},  \ref{CEx}, and \ref{GDx} below, -
and to a certain extent to Section \ref{NMCx}. 
 
When is a category $\E$ suitable for playing the role of  a place, 
or scene, 
where a theory, say axiomatic, of ``spaces'' and their geometry can be developed?

Experience since the 1950s has showed that  many {\em toposes} $\E$ are 
suitable. Thus for instance, the topos of simplicial sets was shown 
to have possibility for being an arena for homotopy theory,  without 
recourse to the real numbers or to the  notion of topological space. 

 A crucial point is that a topos is in many respects like the 
category of sets;  
 in fact, the understanding of ``the'' category of sets  
is distilled out of our experience with categories of spaces, in a broad sense of the 
word. And many 
texts in synthetic differential geometry talk  about the assumed $\E$ {\em as if} it were the category 
of sets, just making sure not to use the law of excluded middle; 
this law holds in the category of {\em abstract} sets (discrete spaces), 
but fails for most other categories of spaces. This is related to the 
contradiction between the discrete and the continuum, see Section 
\ref{Phil} below.

-- So much for the question ``why topos"?''
But why ``ring''? Because we may hopefully use that ring for introducing 
{\em coordinates} and thereby supplement, or even replace, geometric reasoning with algebraic 
calculation. Such ring (a ``number line'')  is however not the central geometric notion in differential 
geometry:

I would like to advance the thesis that the notion of pairs of {\em neighbour} 
points in a space $M$ is a more basic notion in differential geometry;  that 
central differential geometric concepts can be formulated in terms of 
that relation; and that this allows one to present such concept in 
terms of  {\em pictures}, see e.g.\  (\ref{tangx}). There is no notion of ``limiting 
positions'' involved. A ``line type'' ring will be introduced later, 
together with an account of some of the standard synthetic differential 
geometry, and this will also provide ``coordinate'' models for the axiomatics about 
a category of spaces with a neighbour relation, and will thus make this wishful thinking 
 in Chapter I come true.

 The 
guideline for this is the theory developed in Algebraic Geometry by 
K\"{a}hler, Grothendieck, and many others: the first neighbourhood of 
the diagonal of a scheme, or the (first) prolongation space  of a 
smooth manifold  (\cite{KumS} p.\ 52). From these sources, combined 
with standard synthetic differential geometry,   models 
for the axiomatics are drawn (and they are briefly recalled in Section \ref{MAx}).

We use the abbreviation `SDG' for `synthetic differential geometry'. 
  
It is  not a mathematical {\em field}, but a {\em method}. Not 
really a {\em new} one,  it has, as initiated by 
Lawvere,  over several decades by now, 
  contributed, by making the synthetic method more explicit; see also 
Section \ref{Lbx}.

\section{Some differential geometry in terms of the neighbour 
relation}\label{pictx} 

The spaces $M$ considered   in differential geometry come 
equipped with a reflexive symmetric relation $\sim_{M}$, the (first 
order) neighbour relation, often mentioned in the heuristic part of 
classical texts, but rarely made precise, neither how it is defined, nor 
how one reasons with it.

One aim in SDG is to  make precise
how one reasons with the neighbour relation, (and this is done axiomatically); it is not 
a main aim to describe how it is constructed in concrete 
contexts.

The relation $\sim_{M}$ is reflexive and symmetric; it is {\em not} transitive; see the Remark after 
Corollary \ref{31x} why transitivity is incompatible with the 
axiomatics to be presented.

The subobject $M_{(1)}\subseteq M\times M$ defining the relation 
$\sim _{M}$, is 
called the {\em (first) neighbourhood of the diagonal} of $M$. This 
terminology is borrowed from the theory of schemes in algebraic 
geometry, or in differential geometry (see e.g. \cite{KumS}, who call 
it the  {\em  (first) prolongation space} of $M$).

To state some notions of differential-geometric nature, we shall 
talk about the category $\E$ as if it were the category of sets. The 
objects of $\E$, we call ``sets'', or ``spaces''.
If the space $M$ is understood from the context, we write $\sim$ 
instead of $\sim_{M}$. There are also higher order neighbour 
relations $\sim_{2}$, $\sim_{3}, \ldots $ on $M$; they satisfy $(x\sim y) \Rightarrow 
(x\sim_{2}y) \Rightarrow (x\sim_{3}y), \ldots$. The neighbour 
relations 
$\sim_{k}$ will not be transitive, but $x\sim_{k} y$ and $y\sim_{l} z$ will imply 
$x\sim_{k+l}z$. For $x\in M$, we call $\{y\in M\mid y\sim_{k} x\}\subseteq M$ the 
$k$th order {\em monad} around $x$, and we denote it $\M_{k}(x)$. In the 
axiomatics to be presented, it 
represents the notion of $k$-jet at $x$. The first order monad 
$\M_{1}(x)$ will also be denoted $\M(x)$.

The higher neighbour relations will not be 
discussed in the present text, but see e.g. \cite{SGM}.

The spaces one considers live in some category $\E$ of spaces; maps 
in $\E$ preserve the neighbour relations (which is a ``continuity'' 
property). But the neighbour relation on a product space $M \times N$ will be more 
restrictive than the product relation: we will not in general have 
that  
$m\sim m'$ and $n\sim n'$ implies that $(m,n) \sim (m',n')$; see the 
Remark after Proposition \ref{32x}.

For all $x' \in \M (x)$, we thus have by definition $x'\sim x$, and
for sufficiently good spaces, the monad $\M(x)$ will have the 
property that $x$ is the {\em only} point in $\M (x)$ with this 
property.

We present some differential geometric notions that may be expressed in 
terms of the (first order) neighbour relation $\sim$. 
The argument that they comprise the classical 
notions with the corresponding names, may, for most of them,  be 
found in \cite{SGM}.

\subsection{Touching}\label{Touchx}
From the neighbour relation, one derives a fundamental geometric 
notion, namely: what does it mean to say that two subspaces $S$ and $T$ of a 
space $M$ {\em 
touch} each other at a point $x\in S\cap T$?
We take that to mean 
that $ \M_{1}(x)\cap S=  \M_{1}(x)\cap T$.
(The intended interpretation is that $S$ and $T$ are subspaces ``of 
the same dimension''; there is clearly also a notion of, say, 
when a curve touches a surface, which also can be expressed in terms 
of $\sim$.) To ``touch each other at $x$'' is clearly an equivalence 
relation on the set of subspaces of $M$ containing $x$.
(In the intended models, say where $\E$ is a topos containing the 
category of smooth manifolds and $M$ is a smooth manifold, this becomes the relation that $S$ and 
$T$ has {\em first order contact} at $x$.)

Pictures can conveniently be drawn for the 
touching 
notion: in the picture below, $M$ is the plane of the paper, the bullet indicates $x$, 
and the interior of the circle 
indicates $\M_{1} 
(x)$. Ignore the fact that $T$ looks like a line; the notion of {\em 
line} is an  invention of the age of civilization, whereas the notion of {\em 
touching} is known already from pre-civilized stone ages and before. 
So the present section may be thought of as Stone Age Geometry. The 
same applies to Sections \ref{GDx} and \ref{ACx} below.
\begin{equation}
\begin{picture}(100,80)(-18,10)
\put(-12,-2){\qbezier[700](-30,20),(20,10),(70,80)}

\put(6,27){\circle*{3}}
\put(6,27){\circle{15}}
\put(18,20){$\M_{1}(x)$}

\put(-28,10){\line(2,1){90}}
\put(66,80){$S$}
\put(72,54){$T$}

\end{picture}
\label{tangx}\end{equation}

\subsection{Characteristics and envelopes}\label{CEx}

Assume that a space $T$ parametrizes a family $\{S_{t}\mid t\in T\}$  
of subspaces $S_{t}$ of a space $M$. Then for $t_{0}\in 
T$, the  
{\em characteri\-stic} set $C_{t_{0}}$  (at the parameter value $t_{0}$) is the intersection of all the 
neighbouring sets of $S_{t_{0}}$, precisely:
$$C_{t_{0}}:=\cap _{t\sim t_{0}} S_{t},$$
and the {\em envelope} $E$ of the family may be defined as 
\begin{equation}\label{evx}E:= \cup_{t_{0}\in T}C_{t_{0}}= \cup_{t_{0}\in T}\cap _{t\sim 
t_{0}}S_{t}.\end{equation}
Under non-singularity assumptions, $E$ is 
the {\em disjoint} union of the characteristics, i.e.\ there is a 
function $\tau :E\to T$ associating to a point $Q$ of $E$ the parameter 
value $t$ such that $Q\in C_{t}$. So $Q \in \cap_{t'\sim \tau 
(Q)}S_{t'}$.

Let $\tau (Q)=t_{0}$. We would like to prove that {\em $S_{t_{0}}$ {\em touches} $E$ at $Q$}, 
i.e.\ that $ E \cap \M (Q) =  S_{t_{0}}\cap  \M (Q)$ (under a 
``dimension'' assumption on $T$ and the $S_{t}$, commented on below). 
We can in any case prove the inclusion  
$ E\cap \M (Q) \subseteq S_{t_{0}}\cap \M (Q)$. For let $Q'\sim Q$ and $Q'\in 
E$. Then $\tau (Q')\in T$ is defined, and since $\tau $, as any function, preserves 
$\sim$, the assumption  $Q\sim Q'$ implies $\tau (Q')\sim \tau 
(Q)=t_{0}$. Now
$Q' \in \cap_{t\sim \tau (Q')}S_{t}$ by assumption, in particular $Q'\in S_{t_{0}}$.
Therefore  $$ E \cap \M (Q) \subseteq  S_{t_{0}}\cap  \M (Q).$$
(To pass from the proved inclusion to the desired equality would involve a 
dimension argument like: ``the two sets have the same dimension,
so the inclusion implies  equality''. We don't have such 
argument available at this primitive stage; in the Section \ref{NMCx} 
on wave 
fronts, this is 
part of the axiomatics.) 

This leads to an alternative way to describe (but not construct) 
envelopes for such families 
$S_{t}$ of subspaces of $M$. Namely, an 
envelope of such family is {\em a} subspace $E\subseteq M$ such that each 
$S_{t}$ touches $E$, and each point of $E$ is touched by some 
$S_{t}$. This is  also a  classical definition, except that the word 
``touching'' there is defined in terms of differential calculus, not 
available in the Stone Ages. Note the indefinite article ``{\em a} 
subspace''. This ``implicit'' way of describing envelopes is the one 
we use in Section \ref{NMCx} below.

The primary notion in the explicit construction (\ref{evx}) of envelopes is that of {\em characteristic}; the envelope is 
derived from the characteristics. In the literature, 
based on analytic geometry, the characteristic $C_{t_{0}}$ is sometimes, with 
some regret or reservation, defined as ``the limit of the sets 
$S_{t_{0}}\cap S_{t}$ as $t$ tends to $t_{0}$''. 
In \cite{Courant}, Courant (talking about a 1-parameter family of surfaces $S_{t}$ in 3-space, where the 
intersection of any two of them therefore, in non-degenerate cases, is 
a curve) thus writes about a characteristic curve, say $C_{t_{0}}$, for the family: 
``{\em This curve is often referred to in a 
non-rigorous but intuitive way as the intersection of 
 ``neighbouring'' surfaces of the family}'' (p.\ 169) (offering instead: ``{\em If we let $h$ tend 
to zero, the curve of intersection will approach a definite limiting 
position}'' (p.\ 180). What is the topology on the set of subsets which will  justify 
the limit-position notion?)

We shall see (Section \ref{EAx}) that the axiomatics for SDG makes the ``limit'' 
intersection curve rigourous by replacing the dubious {\em limit} with 
the simultaneous intersection of {\em all} neighbouring  surfaces, now 
with ``neighbouring $S_{t}$'' in the strict sense of $t\sim t_{0}$. 
Thus, the ``{\em non-rigourous but} intuitive'' description in Courant's 
text now gets the status: {\em rigourous and} intuitive.

\medskip

\subsection{Wave fronts and rays}\label{NMCx}

Already with the neighbour relation as the only primitive concept, one 
can thus define  the geometric notion (Huygens) of an {\em envelope} of families $S_{t}$ 
of subspaces of a space $M$. Combined with a (weak) notion of 
{\em metric} (distance) on $M$, one can (cf.\ \cite{MSSDG}), by less 
trivial synthetic reasoning 
(and under suitable axioms), recover some of Huygens'  
theory of {\em wave fronts} in geometrical optics: essentially, if 
$B$ is a ``hypersurface'' (in a suitable sense) in $M$, one has an 
envelope $B\, \vdash \, s$ of the family 
of spheres of radius $s>0$ and center on $B$; the Huygens' principle 
states that (for $r$ small enough), this is again a hypersurface, 
``the wave front which $B$ becomes after time lapse $s$''.
 (In particular, Huygens knew that if 
$B$ is a sphere of radius $r$, then $B\, \vdash \, s$ is again a sphere, of 
radius $r+s$.) 

To have a notion of metric, one needs a space of numbers to receive 
the values of the metric, i.e.\ the distances. In the 
intended applications, this will be the strictly positive real 
numbers $R_{>0}$, but only its total strict order $>$ and the 
properties  of the addition operation will be used in the following 
theory; so we are far from being in a situation where a 
coordinatization is used (still, we shall use $R_{>0}$ to denote the assumed 
object that receives the values of the metric). The fact that only strictly positive 
distances are considered means that we  cannot talk about the 
distance from a point to itself; in fact, we cannot talk about the 
distance between a pair of neighbour points. (In the coordinatized 
model of our theory, this has to do with the fact that the square 
root function is not smooth at $0$.) When we say that two points are 
{\em distinct}, we thus imply that their distance is defined (hence 
positive).

With a metric on $M$, we can define {\em spheres}: if $a\in M$ and 
$r\in R_{>0}$, the sphere $S(a,r)$ with center $a$ and radius $r>0$ 
is the set $\{b\in M \mid ab =r\}$, where ''$ab$'' is short notation 
for the distance between $a$ and $b$. So $ab=ba$. We assume that, as in Euclidean 
geometry, the center $a$ and the radius $r$ can be reconstructed from 
the point set $S(a,r)$. No triangle inequality is used in the 
following.

Combining the two primitive notions: neighbours and metric, we can 
then define  the notion of  contact element $P$: A {\em contact element} at 
$b\in M$ is a subset of the form $P=\M (b)\cap B$, where $B$ is a 
sphere with $b\in B$. The same contact element may be presented in 
$ \M (b)\cap B'$ for many other spheres $B'$, but all these 
spheres touch each other at $b$, since $\M (b) \cap B = P = \M(b) 
\cap B'$.

Since a contact element at $b$ has $P \subseteq \M (b)$, one has 
that $b$ is 
neighbour of all the points in $P$.  We assume that $b$ 
is the {\em }only 
point in $P$ with this property. So $b$ can be reconstructed from the 
point set $P$; we may all it the {\em focus} of $P$, to avoid saying 
``center''. 

In the intended application, where $M$ is a smooth manifold, the set 
of contact elements make up the projectivized cotangent bundle of 
$M$.

Note that in the classical theory, any contact element $P$ at $b$, say $\M(b) 
\cap B$ (where $b\in B$), is a one-point set, $P=\{b\}$ , since $\M(b)$ is so; whereas with a 
non-trivial $\sim$, there is much more information in $P$: it 
generates a non-trivial perpendicularity relation. Namely, for $c$ 
distinct from $b$ (and (hence) from any $b'\sim b$), we say that $c$ 
is {\em perpendicular to $P$}, or 
$c\perp P$, if for all $c' \in P$, we have $bc'=bc$ (where $b$ is the 
focus of $P$). (For a  
trivial $\sim$, all points distinct from $b$ are perpendicular to 
$P=\{b\}$.) 

There are two basic structures in geometrical optics, (light-) {\em 
rays} and {\em wave fronts}. These can be described in the present 
framework. The {\em rays} in $M$ are certain (open) half lines, parametrized 
 by $R_{>0}$, described more precisely below in terms of a 
 collinearity condition. Wave fronts here occur in the present 
context  as {\em (hyper-)surfaces}; the rude notion of {\em hypersurface} we are considering is the 
following: it is  a subset of $M$ 
which ``made up of 
contact elements'', i.e.\ it is a subset $B 
\subseteq M$ such that for each $b\in B$, the set $\M (b)\cap B$ is a contact 
element  (necessarily with focus $b$). So in particular, a sphere is a hypersurface.

What makes synthetic reasoning about rays and wave fronts possible, 
is an analysis about how spheres may touch. In the deductions in 
\cite{MSSDG}, this analysis takes form of two axioms, one for 
``external'' touching and one for ``internal'' touching. We state 
them below, noting that they are refinements of theorems of Euclidean 
geometry (in the sense that ``{\em touch}'' has a refined meaning). Thus for external touching: 

{\em Two spheres touch (externally) if 
the sum their radii is the distance between their centers}. 

Here is 
the picture for external touching; the spheres are $A=S(a,r)$, 
$C=S(c,s)$, they touch at $b$. The two other dots represent the 
centers $a$ and $c$:

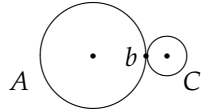
$$
\begin{picture}(100,40)(0,10)
\put(-2,17){$A$}
\put(30,30){\circle{40}}
\put(30,30){\circle*{2}}

\put(64,17){$C$}
\put(58,30){\circle{15}}
\put(58,30){\circle*{2}}
\put(50,30){\circle*{2}}
\put(42,26){$b$}

\end{picture}
$$
With the refined touching notion derived from $\sim$, here is how 
this basic fact gets formulated: given $A=S(a,r)$ and $C=S(c,s)$. If $r+s 
= ac$, then the spheres $A$ and $C$ touch at a unique 
point $b$, and this $b$ is characterized by 
\begin{equation}\label{A11x}\mbox{ $ab+bc=ac$; and for all } b'\sim b,  \mbox{  we have   } (ab'=ab) 
\Leftrightarrow 
 (b'c =bc).\end{equation}
 This is essentially the basic Axiom (together with a similar axiom 
 for internal 
 touching), except that we weaken it by replacing the 
 $\Leftrightarrow$ in (\ref{A11x}) by  $\Rightarrow$: 
 \begin{equation}\label{A12x}\mbox{ $ab+bc=ac$; and for all } b'\sim b,  \mbox{  we have   } (ab'=ab) 
\Rightarrow 
 (b'c =bc).\end{equation}
 This replacement is, in the intended models, justified by a dimension argument, as alluded to 
 in Section \ref{CEx}.  Note that (\ref{A12x}) can also be expressed: 
 $c$ is a characteristic point (in the sense of Section \ref{CEx}), for parameter value $b$, of the family $S(b',s)$, as $b'$ ranges 
 over $P= \M(b)\cap  S(a,r)$, or, as $b'$ ranges over $S(a,r)$. For, $b'c=s$ is equivalent to $c\in 
 S(b',s)$.
 
We give an equivalent formulation of (\ref{A11x}),  and also of the 
corresponding way of writing the Axiom for internal touching; $a$, 
$b$ , and $c$ denote points in $M$, and $s$ denote an  element $\in R_{>0}$. Note that for 
internal touching, there is no restriction on $s\in R_{>0}$.

{\em Given $a$, $c$, and $s$, with $s< ac$. Then $\exists ! b$ such that $S(a,ac-s)$ and 
$S(c,s)$ touch at $b$.
Given $a$, $b$, and $s$. Then $\exists ! c$ such that $S(a,ab+c)$ and 
$S(b,s)$ touch at $c$.}

Since touching of two spheres is either internal or external, it is straightforward 
that one can (transversally) 
{\em orient}   a contact element $P$ in two ways, and one then can  divide the 
class of points perpendicular to $P$ in two classes, those on the 
``outer'' side and those on the ``inner side''. They are the two 
{\em rays} defined by $P$: given an orientation of $P$ and an $s\in 
R_{>0}$, we let 
$P\, \vdash \,  s$ denote the unique point on the outer side perpendicular 
to $P$ and whose distance 
to $b$ is $s$ (where $b$ denotes the focus of $P$). 
To jusitify the word ``ray'', note that 
the ray generated by $P$ is a point set, bijectively parametrized by 
$R_{>0}$; and furthermore, any three distinct points (taken in 
suitable order) on this ray are {\em 
collinear}:

Collinearity is a notion 
which, when $\sim$ is trivial, may be formulated purely 
in terms of the metric, and it forms the basis of Busemann's theory of 
geodesics, cf.\  \cite{Busemann}. Three points $a,b,c$ (say, distinct) 
are  classically and in loc.\ cit.\ called collinear if $ab+bc=ac$.
The stronger collinearity property which applies to the rays in the 
present theory is that $ab+bc=ac$ {\em and} that $S(a,ab)$ touches 
$S(c,bc)$ in $b$; equivalently, if (\ref{A11x}) (or (\ref{A12x})) 
holds (with $r=ab$, $s=bc$).

To give a hypersurface an {\em orientation} is to give each of its contact 
elements an orientation.
Let $B$ be an oriented hypersurface, and let $s>0$. Let us denote by 
$B \, \vdash \, s$ the set of points of the form $\M (b)\, \vdash \, s$, i.e.\ 
the envelope (in the explicit sense of (\ref{evx})) of the spheres $S(b,s)$ as $b$ ranges over $B$. Since 
the distance of $b$ and $\M (b)\, \vdash \,  s$ is $s$, we would like to 
think of $B\, \vdash \, s$ as the parallel hypersurface to $B$ at distance 
$s$; however, it may not be a hypersurface, as is well known in geometry, 
even when $M$ is the Euclidean plane: there may be 
self-intersections, cusps, etc.\ if $B$ is concave. But unless $B$ is 
very crinkled, one will for sufficienty small $s$ have that the map 
$s\mapsto \M(b)\, \vdash \, s$ is a bijection $B \to B\, \vdash \, s$. The 
version of Huygens' principle we can prove synthetically 
(cf.\ \cite{MSSDG}) is:  

{\em Assume that $B$ is an oriented  hypersurface, and that $s>0$ is so 
that the map $B \to B \,\vdash \,  s$ described is a bijection. Then 
$B\, \vdash\, s$ is again a hypersurface.}

If $B=S(a,r)$ is a sphere, $B\, \vdash \, s$ will be the sphere 
$B=S(a,r+s)$ for one orientation of $B$, and will, for the other 
orientation,   be the 
sphere $S(a,r-s)$ (provided $s<r$).

\subsection{Geometric distributions}\label{GDx}

A (geometric) {\em distribution} on $M$ is a reflexive symmetric 
relation $\approx$ refining 
$\sim$ (i.e.\ $x\approx y$ implies $x\sim y$). It is called {\em involutive} if it satisfies, 
for all $x,y,z $ in $M$:
$$(x\approx y) \wedge (x\approx z) \wedge (y\sim z) \mbox {\quad  
implies \quad } y\approx z.$$ 
A relevant picture is the following; single lines indicate the 
neighbour relation $\sim$, double lines indicate the assumed ``strong'' 
neighbour relation $\approx$.

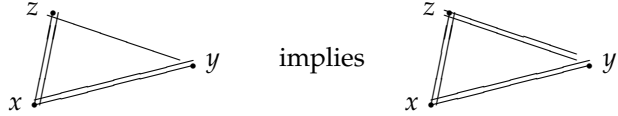
\begin{equation*}
\begin{picture}(200,75)(0,-10)
\put(20,6){\line(4,1){60}}
\put(20,8){\line(4,1){60}}
\put(20,6){\line(1,5){7}}
\put(22,6){\line(1,5){7}}
\put(10,4){$x$}
\put(85,21){$y$}
\put(17,40){$z$}
\put(20,6){\circle*{2}}
\put(80,21){\circle*{2}}
\put(27,41){\circle*{2}}
\put(25,40){\line(3,-1){50}}
\put(110 ,21){\mbox{ implies }} 
\put(170,6){\line(4,1){60}}
\put(170,8){\line(4,1){60}}
\put(170,6){\line(1,5){7}}
\put(172,6){\line(1,5){7}}
\put(160,4){$x$}
\put(235,21){$y$}
\put(167,40){$z$}
\put(170,6){\circle*{2}}
\put(230,21){\circle*{2}}
\put(177,41){\circle*{2}}
\put(175,40){\line(3,-1){50}}
\put(175,42){\line(3,-1){50}}
\end{picture}
\label{simplexc}
\end{equation*}

For instance, if $f:M\to N$ is any map between spaces, the relation 
on $\approx$ on  
$M$ 
defined by $x\approx y$  iff  $(x\sim y)\wedge (f(x)=f(y))$ is a distribution, 
in fact an involutive one.

\medskip
An {\em integral subset} of a distribution $\approx$ is a subset 
$F\subseteq M$  such that on $F$, the relations $\sim$ and 
$\approx$ agree.   An important 
integration theorem in differential geometry is Frobenius' Theorem, 
whose conclusion is that for an {\em involutive} distribution, there 
exist maximal connected
integral subsets (leaves). 

Such integration results can usually not be {\em proved} in the 
context of SDG (even the 
very formulation may require some further primitive concepts), since they in a more serious way depend on limits and 
on completeness of the real number system. Sometimes, SDG can {\em 
reduce} one integration result to another; this is also an old 
endeavour in classical differential geometry, e.g.\ Lie has many 
results about which differential equations can be solve {\em by 
quadrature}, i.e.\ by reduction to existence of anti-derivatives.

\medskip

\noindent {\bf Example.} The following is meant as a sketch of an 
(involutive) distribution in the plane. Consider  

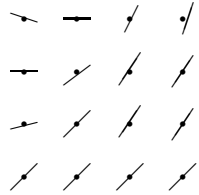
\begin{picture}(150,110)(-96,20)

\put(40,40){\circle*{2}}
\put(35,35){\line(1,1){10}}

\put(60,40){\circle*{2}}
\put(55,35){\line(1,1){10}}

\put(80,40){\circle*{2}}
\put(75,35){\line(1,1){10}}

\put(100,40){\circle*{2}}
\put(95,35){\line(1,1){10}}

\put(40,60) {\circle*{2}}
\put(35,58){\line(4,1){10}}

\put(60,60) {\circle*{2}}
\put(55,55){\line(1,1){10}}

\put(80,60){\circle*{2}}
\put(76,55){\line(2,3){8}}

\put(100,60){\circle*{2}}
\put(96,54){\line(2,3){8}}

\put(40,80) {\circle*{2}}
\put(35,80){\line(1,0){10}}

\put(60,80) {\circle*{2}}
\put(55,75){\line(4,3){10}}

\put(80,80){\circle*{2}}
\put(76,75){\line(2,3){8}}

\put(100,80){\circle*{2}}
\put(96,74){\line(2,3){8}}

\put(40,100){\circle*{2}}
\put(35,102){\line(3,-1){10}}

\put(60,100){\circle*{2}}
\put(55,100){\line(1,0){10}}

\put(80,100){\circle*{2}}
\put(78,95){\line(1,2){5}}

\put(100,100){\circle*{2}}
\put(100,94){\line(1,3){4}}

\end{picture}

\noindent In this picture, the ``line segments'' are the   
$\approx$-monads $\M_{\approx}(x):= \{y\mid y\approx x\}$  
 around (some of) the points $x$ (drawn as dots) of $M$. But note that
 the notion of ``line'' has not yet entered in our 
vocabulary, let alone coordinate systems like $R \times R$; when such 
things are present, an ordinary first order differential equation
$$y'=F(x,y),$$ as in the Calculus Books, gives rise to such a picture, 
known as the ``direction field'' of the equation: through each point $(x,y)\in R\times R$, one draws a ``little'' line 
segment $S(x,y)$ with slope $F(x,y)$.

The ``integral subsets'' of a distribution of this kind are 
essentially (the graphs of) the solutions of the differential 
equation.

I cannot draw a good picture of a non-involutive distribution: 
paper is  2-dimensional.  But in three dimensions: consider the 
scales of a ripe pine cone, and extrapolate radially.

If $M$ carries a metric (in the sense of Section \ref{NMCx}), it makes 
sense to say that a distribution is ``of codimension 1'' if all the 
$\approx$-monads are contact elements. The two specific examples 
mentioned have this property.

\subsection{Affine connections}\label{ACx}
An {\em affine connection} on  a space $M$ is a law $\lambda$ which completes 
any  configuration $(x,y,z)$ consisting of three points $x$, $y$, and 
$z$ with $x\sim y$, $x\sim z$ by a fourth point $\lambda 
(x,y,z)$ with $y\sim \lambda (x,y,z)$ and $z\sim \lambda(x,y,z)$:

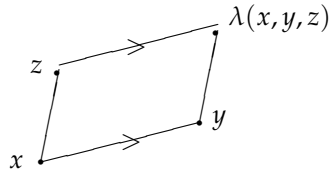
\begin{equation}
\begin{picture}(140,75)(0,-10)
\put(20,6){\line(4,1){60}}
\put(20,6){\line(1,5){7}}
\put(50,11){$>$}
\put(52,47){$>$} 

\put(80,21){\line(1,5){7}}

\put(27,43){\line(4,1){60}}
\put(8,4){$x$}
\put(85,21){$y$}
\put(16,40){$z$}

\put(80,21){\circle*{2}}
\put(20,6){\circle*{2}}

\put(26,40){\circle*{2}}
\put(86,55){\circle*{2}}
\put(91,58){$\lambda(x,y,z)$}

\end{picture}
\label{parallelogram}\end{equation}
expressing  ``infinitesimal parallel transport of $z$ along 
$\overline{xy}$'', or ``constructing an 
infinitesimal parallelogram''. (We assume $\lambda(x,x,z)=z$, and 
$\lambda(x,y,x)=y$).) The connecting lines 
indicate the assumed neighbour relations. We us different 
signature for the edges $xy$ and $xz$, since we do not assume the 
symmetry condition $\lambda(x,y,z)= \lambda (x,z,y)$. If symmetry 
holds, $\lambda$ 
is  called a {\em symmetric} or {\em torsion free} connection.

\medskip

A {\em geodesic} for a given torsion free affine connection on $M$ is a subset 
$S\subseteq M$ which is stable under $\lambda$ in the sense that if 
$x\sim y$ and $x\sim z$ with $x$, $y$ and $z$  in $S$,  then
$\lambda (x,y,z) \in S$.

\medskip

The {\em curvature} of  an affine connection may be described 
combinatorially by asking the question: what happens if we transport 
$z\sim x_{0}$ around a  circuit from $x_{0}$ to $x_{1}$, then from 
$x_{1}$ to $x_{2}$, and finally from $x_{2}$ back to $x_{0}$? This 
makes sense whenever $x_{0}\sim x_{1}\sim x_{2}$ {\em and} $x_{0}\sim 
x_{2}$ (the latter requirement  is not automatic: the relation $\sim$ is not 
transitive). The result of such circuit transport gives a new point 
$z'\sim x_{0}$; thus the ``infinitesimal 2-simplex'' $(x_{0},x_{1},x_{2})$ provides an 
automorphism  $z\mapsto z'$, denoted $R(x_{0},x_{1},x_{2})$, of the pointed set $\M (x_{0})$; this is the  
curvature of $\lambda$, more precisely, the curvature of $\lambda$ is the law which to an 
infinitesimal 2-simplex $(x_{0},x_{1},x_{2})$ associates the 
described automorphism of $\M (x_{0})$. If  this automorphism is the 
identity map for all infinitesimal 2-simplices, the connection is 
called {\em flat}.
(Any affine  connection on a 1-dimensional space $M$ is flat. 
One may even experiment with this as a definition of ``$M$ is of 
dimension (at most) 1''.) 

\subsection{Differential forms}\label{DFx} Differential forms are, in analytic 
differential geometry, certain functions taking values in a {\em ring} $R$ of quantities (or in a 
{\em module} over $R$), but are  in the present context (equivalent to) a 
special case of
 a more primitive, non-quantitative, kind of thing: Thus, in SDG, one 
may, for any group $G$, define ``(combinatorial) $G$-valued $k$-form on a 
space $M$'' to mean a ``function  $\omega$, which takes as  input  
infinitesimal $k$-simplices ($k+1$-tuples of 
mutual neighbour points in $M$), and returns as output  elements in 
$G$''. One imposes the normalization condition that $\omega (x_{0}, 
\ldots ,x_{k})=e$ whenever two of the $x_{i}$s are equal (where $e$ 
denotes the neutral element of $G$). A $G$-valued $0$-form on $M$ is 
then just a function $f:M \to G$; it has a ``coboundary'' $df$, which 
is a $G$-valued 1-form, defined by $df(x_{0},x_{1}):= 
f(x_{0})^{-1}\cdot f(x_{1}) 
$. A $G$-valued 1-form $\omega$ on $M$ has a coboundary 
$d\omega$, which is a $G$-valued 2-form defined by $$d\omega (x_{0},x_{1}, x_{2}):= \omega(x_{0}, 
x_{1})\cdot \omega (x_{1},x_{2})\cdot \omega (x_{2},x_{0}).$$
The 1-form $\omega$ is {\em closed} if $d\omega$ is constant $e$. The 
1-form $df$ is always closed.

The group $G$ carries a canonical closed $G$-valued 1-form, namely 
$df$, where $f:G \to G$ is the identity function. This is the 
Maurer-Cartan form of $G$. 

If there is given data identifying all the $\M (x)$ of a 
given manifold $M$ with each other, then the curvature of an affine connection $\lambda$ 
on $M$ may be seen as a 2-form with values in the automorphism group of 
the pointed set
$\M (x_{0})$ (for some, hence any, $x_{0}\in M$). (Alternatively, one 
gets a 2-form ``with local coefficients''; then no identification 
data is needed.)

For most $G$, we have that $G$-valued differential forms are {\em alternating}: 
interchanging two of the input entries implies inversion of the value 
of the form.  In 
particular  $$\omega (x_{0},x_{1})^{-1}= \omega (x_{1}, x_{0}).$$ Such 
1-form $\omega $ on $M$ then defines a geometric distribution on $M$ 
by saying $x\approx y$ iff $\omega (x,y)=e$. If $\omega$ is closed, 
the $\approx$ which is defined by $\omega$ is involutive.

There is a relationship between combinatorial group valued 1-forms, on the one 
hand, and the general notion of connection in a fibre bundle,  or  in a groupoid, 
on the other. This we expound in Section \ref{Connx} below.

In case the value group $G$ is commutative (additively written), there 
are, for good $M$ and $G$, a natural bijection between combinatorial 
$G$-valued 
forms, and the standard multilinear alternating forms on $T(M)$, the 
tangent bundle of $M$, see \cite{SDG} I.18.

\section{Neighbours in the context of Euclidean geometry}\label{NCEGx}
In this Section, we move from the Stone Age into  the era of Civilization, and
assume that some classical Euclidean geometry (plane, say) is available in 
a space $E$ (with a given neighbour relation $\sim$). In particular, 
there are given subsets called points and lines; they are  {\em affine} 
subspaces of $E$ (without yet assuming the existence of a ``number'' 
line $R\subseteq  
E$, i.e.\ a line equipped with  a commutative 
ring structure). 

Then we can be more explicit about our wishes for the compatibilities 
between the Euclidean notions and the combinatorics of the neighbour 
relation. We refrain from calling these wishes for 
``Axioms'', since they (for the coordinate spaces $R^{n}$ built on 
$R$) lead to and) are subsumed in a more complete 
comprehensive axiom scheme later on; so we call these wishes for 
``Principles''.

There are also some incompatibilities, essentially because in 
Euclid, the law of excluded middle is explicitly used. Thus, in 
Euclid, a curve, say a circle, has {\em exactly one} point in common with 
any of its tangents, so that the picture (\ref{tangx}) (with $S$ as 
part of a circle) is an illusion 
for Euclid; already the contemporary  Greek philosopher Protagoras is said to 
have ridiculed Euclidian geometry for insisting on the ``only one 
point''-idea, which seemed to him  to go against experience. 
In Euclid's geometry, $\M (x)$ is always just the  one-point set 
$\{x\}$.
Certainly, the following principle is incompatible with such a small 
$\M (x)$; in the terminology of Chapter \ref{pictx}, this says that 
two lines which touch each other at some point are equal.
 
\medskip

\noindent {\bf Principle.}  {\em Given two lines $l_{1}$ and $l_{2}$ 
in a plane $E$. 
Let $x\in  l_{1} \cap l_{2}$. Then $$\M (x) \cap l_{1} = 
\M(x)\cap l_{2}\mbox{ implies } l_{1} = l_{2}.$$}
A subspace $C\subset E$ is called a {\em curve}  if for each $x\in C$, 
there exists a line $l$ such that $ \M (x)\cap l =  \M (x)\cap C$, 
i.e.\ a line which touches $C$ at $x$; 
such  a line is  unique, by the Principle. This line then deserves the 
name: the {\em tangent} of $C$ in $x\in C$.
In the picture (\ref{tangx}), if $T$ is a line (as the picture 
suggests), then the picture says that this line is the tangent to the curve $S$ at $x$.

For any curve $C$, the family of its tangents $T_{x}$ ($x\in C$) is a parametrized 
family, parametrized by the points of curve $C$.

\begin{prop}Any curve $C$ is contained in the  envelope of its family of tangents.
\end{prop}
{\bf Proof.} For $z\in C$, let $T_{z}$ denote the tangent  $C$ at $z$.
Let $x\sim y$ be points in $C$. So $y\in \M (x)\cap C = \M (x) \cap 
T_{x}$; so $y\in T_{x}$. Similarly, $x\in T_{y}$. So 
$$x\in \bigcap _{y\in \M (x)\cap C}T_{y},$$
which is to say that $x$ belongs to the characteristic set (for 
parameter value $x$) of the family of tangents. Hence it belongs to 
the envelope of the family.

\medskip

Let $M$ and $N$ be spaces (objects in $\E$). 
It will not in general be the case that $x\sim x'$ and  $y\sim y'$ in 
$N$ implies $(x,y) \sim (x',y')$ in $M\times N$ (although the 
converse implication will hold, since the projections, like any other 
map, preserve the assumed neighbour relations $\sim$). But if
$f:M\to N$ is a map, 
then we also have  that 
\begin{equation}\label{graphx}x\sim x' \mbox{ in $M$  iff  } 
(x,f(x))\sim (x',f(x'))\mbox{ in }M\times N ;\end{equation}  
for, the map $M \to M\times N$ given by $x\mapsto (x,f(x))$ preserves, 
like any map, the neighbour relation.

In the classical geometry of conics, consider a parabola. Then the tangent 
line at the apex  is perpendicular to the axis of the parabola:


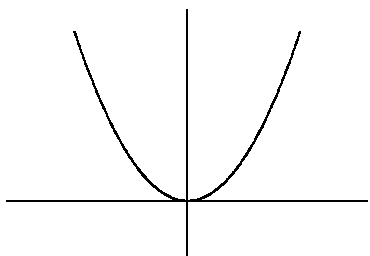
\setlength{\unitlength}{.3cm}

\begin{picture}(100,12.5)(-5,-15)
\put(12,-10){\qbezier[500](-5,5),(0,-10),(5,5)}
\put(12,-12.5){\line(1,0){8}}
\put(12,-12.5){\line(-1,0){8}}
\put(12,-15){\line(0,1){11}}

\end{picture}
When coordinates are introduced in the plane, by making ``the'' geometric 
line $R$ into a commutative ring, we may consider in particular the 
parabola $P$ given as the graph of $y=x^{2}$. The axis of $P$ is the 
$y$-axis, and the tangent line at 
the apex of $P$ is the $x$-axis $X$.  It follows from (\ref{graphx}) that 
$x\sim 0$ implies 
$(x,x^{2})\sim  (0,0)$. Since  $\M (0,0)\cap P = \M (0,0) 
\cap X$, we conclude that $(x,x^{2}) \in X$, which implies that 
$x^{2}=0$. Thus
 \begin{equation}\label{protx}x\sim 0\mbox{ implies }x^{2}=0,
\end{equation}
or, writing $D\subseteq R$ for $\{x\in R \mid x^{2}=0\}$, this says 
that $\M (0)\subseteq D$. (The other inclusion will be our definition 
of $\sim$ in this coordinate model.)

Next consider $(x,y)\sim (0,0) \in R^{2}$. Since the projections 
$R^{2}\to R$ preserve $\sim$, we conclude $x\sim 0$ and $y\sim 0$, so 
$(x,y)\in D \times D$, i.e.\ $x^{2}=y^{2}=0$; but we can say more, 
namely that $x\cdot y =0$. For, the addition map $R\times R \to R$ 
preserves, like any map, the neighbour relation, so $(x,y) \sim 
(0,0)$ implies $(x+y)^{2}=0$. But $(x+y)^{2}=x^{2}+y^{2}+2x\cdot y$. 
The two first terms we already know are $0$, hence so is $2x\cdot y$, 
and since $2$ is invertible, we conclude $x\cdot y =0$; thus
\begin{equation}\label{protxx} (x,y)\sim (0,0) \mbox{ implies }
x^{2}=y^{2}= x\cdot y =0.
\end{equation}

We embark in the following Section on a more serious investigation of how synthetic notions 
like $\sim$ can be conveniently coordinatized, by suitable 
axiomatization of properties of the ring $R$.

\section{Coordinate geometry, and the axiomatics}\label{CGx}

It is not the intention of SDG to avoid using the wonderful tool of 
coordinates. So we now embark on  the interplay between an assumed 
neighbour relation on the spaces, and an assumed 
basic geometric line $R$ with a commutative ring structure.

The reason we did not start there, is to stress that the 
``arithmetization'' in terms of $R$ is a {\em tool}, not the {\em subject matter}, 
of geometry. This also applies in differential geometry, which  has 
some important aspects without any $R$ (as illustrated by the material in 
Chapter \ref{pictx} and partly in Chapter \ref{NCEGx}); so in particular, it has a life 
without the ring ${\mathbb R}$ of real numbers, who sometimes thinks 
of himself as being the owner and boss of the company. 

The scene of SDG in its present form is thus a category $\E$ (whose object we call {\em spaces} or 
{\em sets}), together with a commutative ring object $R$ in it.  But 
$\E$ is not the category of {\em discrete} sets, so some of the 
logical laws valid for the category of discrete sets, like the law of excluded 
middle, cannot be used.
In differential geometry, whose maps are {\em smooth} maps, the law 
of excluded middle does anyway not apply; it would immediately lead out of the 
smooth world, like when one attempts to construct the absolute value 
function $x\mapsto |x|$ on the number line. 

Nevertheless, we shall talk about the objects and maps of $\E$ as if 
they were sets; just recall that they are not {\em discrete} 
sets\footnote{See the discussion in Section \ref{Phil}}. This is 
a basic technique in modern mathematics, more or less explicitly 
used in many other contexts. We shall not say more about it here. Basic concepts for making 
the technique explicit are Cartesian closed categories, or even 
better, locally 
Cartesian closed categories, in particular toposes  (when talking about 
``families'' of objects, as in the 
discussion above on envelopes). (There is some explicit description 
of the technique relevant for SDG  in Part II of \cite{SDG} and in 
Appendix A2 in \cite{SGM}.) 

\medskip

 The axioms concern $R$; the category $\E$ should just have 
sufficiently good properties. The 
maximal thing wanted is that $\E$ is a topos, but less will often do. Thus, to get hold of an object like 
the unit circle
$\{(x,y)\in R^{2}\mid x^{2}+y^{2}=1\}$, one needs only that $\E$ has finite 
limits; the circle then is a subobject of $R\times R$ given as the 
equalizer of two particular maps $R\times R \to R$. (In fact, the term 
``equalizer'' came from such equational conditions as 
$x^{2}+y^{2}=1$.)

For simplicity, we therefore in the following assume that $\E$ is a 
topos; and that $R$ is a commutative ${\mathbb 
Q}$-algebra in it. The intuition and terminology is: $R$ is the number 
line; and also: $R$ is the ring of scalars.

\subsection{The axiomatics}

The axiom for such data, which is at the basis of the form of SDG 
considered here, is an axiom{\em -scheme}\footnote{often referred to 
as the general KL axiom, for ``Kock-Lawvere'', cf.\ e.g.\ \cite{MR} 
or \cite{Lav}}, with one axiom for each Weil 
algebra; a Weil algebra is a finite dimensional 
commutative  algebra (over ${\mathbb Q}$, for the present purpose), 
where the nilpotent elements form 
an ideal of codimension $1$. The name `Weil algebra' is used because they 
were introduced in the ``Points proches''-paper by A.\ Weil, 
\cite{Weil}, whose 
aim was related to the one we present here. The simplest non-trivial 
Weil algebra is the 
``ring of dual numbers'' ${\mathbb Q}[\epsilon ] = {\mathbb Q}[X]/(X^{2})$, which is 
2-dimensional over ${\mathbb Q}$.

Concerning $R$, we have already seen in (\ref{protx}) that $\M(0)\subseteq \{x\in R 
\mid x^{2}=0\}$. The latter object we call $D$, as at the end of 
Chapter \ref{NCEGx}. To relate the combinatorics of $\sim$ with the 
algebra of $R$, we postulate the converse inclusion $D\subseteq \M (0)$. 
It then follows that $x\sim y$ in $R$ iff $(y-x)^{2}=0$.

 The simplest instantiation of the 
axiom scheme concerns $D$. It can be seen as the instantiation of the 
axiom scheme for the two-dimensional Weil algebra 
 $${\mathbb Q}[\epsilon ]:= {\mathbb Q}[X]/(x^{2}).$$

\medskip

\noindent {\bf Axiom 1.} {\em  Every map $f:D \to R$ is of the form 
$d\mapsto a+d\cdot b$
for unique $a$ and $b$ in $R$.}

\medskip
\noindent This has to be true  {\em with parameters}, thus if $f:I\times D \to R$ 
is an $I$-parametrized family of maps $D \to R$, then the $a$ and $b$ 
asserted by the axiom are likewise 
$I$-parametrized points of $R$, i.e.\ they are maps $I \to R$. In a 
Cartesian closed category $\E$, the ``true with parameters'' 
 follows from a more succinct property, namely
the property that the map $R\times R \to R^{D}$, given by $(a,b) \mapsto [d\mapsto 
a+d\cdot b]$, is invertible. Thus, the axiomatics for SDG is simpler to 
state under the assumption that the category $\E$ is Cartesian 
closed (although the idea and logic of parametrized families can also 
be made precise, even without Cartesian closedness).

 Cartesian closedness of $\E$ is an aspect of talking about 
the objects of $\E$ as if they were  sets.

We leave to the reader to prove (using $D$  as a space of 
parameters)
\begin{cor}\label{31x}Every map $f: D\times D \to R$ is  of the form
$(d_{1},d_{2}) \mapsto a + d_{1}\cdot b_{1} + d_{2}\cdot b_{2}  + 
 d_{1}\cdot d_{2}\cdot c$, for unique $a$, $b_{1}$, $b_{2}$, and $c$ in $R$.
\end{cor}
In rough terms, since $R^{D}\cong R^{2}$, it follows that 
$(R^{D})^{D}\cong (R^{2})^{D}\cong (R^{D})^{2}\cong (R^{2})^{2}\cong 
R^{4}$. In itself, the Corollary also appear as an instantiation of the axiom 
scheme, namely for the four-dimensional Weil algebra $${\mathbb 
Q}[\epsilon_{1},\epsilon_{2}] := {\mathbb Q}[X_{1}, X_{2}]/(X_{1}^{2}, 
X_{2}^{2}).$$

\medskip

\noindent {\bf Remark.} The relation $\sim$ defined in terms of $D$ 
cannot be transitive. For, transitivity is easily seen to be equivalent to $D$ 
being stable under addition, and hence (using that $2$ is invertible) 
that $d_{1}\in D$ and $d_{2}\in D$ implies $d_{1}\cdot d_{2}=0$. But 
this contradicts the uniqueness of the coefficient $c$ in the above 
Corollary. So Axiom 1 implies that $\sim $ is not transitive.

\medskip 

We shall not be explicit how one goes from a (finite presentation 
of) a Weil algebra to the corresponding Axiom (see \cite{SDG} I.16). The 
reader may guess the pattern from the examples given.

 From the uniqueness assertion in Axiom 1 one derives

\medskip

\noindent{\bf Principle of cancelling universally quantified $d$'s}: 
let $r, s \in R$. Then:
\newline
\indent {\em If $d\cdot r = d\cdot s$ for all $d\in D$, then $r=s$}.

\medskip 

\noindent In the classical treatment, any individual $x\neq 0$ in $R$ 
is cancellable, i.e.\ it 
has the property  that it 
detects equality; $x\cdot r =x\cdot s$ implies $r=s$; for, in the 
classical treatment, $R$ is a field, so $x\neq 0$ implies that $x$ is 
invertible. On the other hand, in SDG, no individual $d\in D$ can be 
cancellable; for, any such $d$ is nilpotent. This, for some intuition,  
means that $d$ is very small, ``infinitesimal''. So none of these small 
elements individually have   the strength that they can detect equality; 
but when the small elements join hands, they can. Collective strength, 
of all the small together,  replaces the strength of any individual.

\medskip

Another consequence of the Axiom 1 is that the beginnings of 
differential {\em calculus} become available: given $f:R \to R$, one 
applies, for each $x\in R$, the axiom to the function $d\mapsto 
f(x+d)$; so one gets for each $x$ that there are unique $a$ and $b$ 
such that  $f(x+d)=a+d\cdot b$ for all $d\in D$. The $a$ and $b$ depend on the $x$ 
chosen, so write them $a(x)$ and $b(x)$, respectively. By setting 
$d=0$, we conclude $a(x)=f(x)$; but $b(x)$  deserves a new name, we 
call it $f'(x)$, so for all $d\in D$, we have the exact  ``Taylor expansion''
\begin{equation}\label{taylor1x}f(x+d)= f(x)+d\cdot f'(x)\mbox{ for 
all } d \mbox{ with } d^{2}=0.\end{equation}
And this property characterizes $f'(x)$, by the principle of 
cancelling universally quantified $d$s.

Since such Taylor expansion holds also with parameters, one 
also gets partial derivatives for functions in several variables, by  
considering the variables, except one, as parameters. See 
(\ref{partialx}) below for an example.

\medskip

\noindent{\bf Remark.} For differential {\em calculus}, there are 
other synthetic/axiomatic theories available: e.g.\ the 
``Fermat''-axiom (suggested by Reyes), see e.g.\  \cite{MR} VII.2.3; 
the axiomatics of ``differential categories'' (cf.\ \cite{BCS}, 
\cite{CC}, and references therein); and  the 
``topological differential calculus'' (cf.\ \cite{Be}, and references therein).

\medskip

The Corollary \ref{31x} could be seen as an instantiation of the 
general 
axiom scheme; a more interesting instantiation  of the axiom scheme comes about by 
considering the  three-dimensional Weil algebra  
$${\mathbb Q}[\epsilon_{1},\epsilon_{2}]/(\epsilon_{1}\cdot \epsilon 
_{2}):={\mathbb Q}[X_{1},X_{2}]/(X_{1}^{2},X_{2}^{2}, 
X_{1}\cdot X_{2}). $$
To state the Axiom, let $D(2) \subseteq R^{2}$ be given as 
$$\{(d_{1},d_{2}) \in R^{2}\mid d_{1}^{2}= d_{2}^{2} = d_{1}\cdot 
d_{2}=0\}.$$
 (Clearly, $D(2) \subseteq D\times D$. Note that $D(2)$ is defined by 
the equations 
occurring in (\ref{protxx}).) Then

\medskip

\noindent {\bf Axiom 2.} {\em  Every map $f:D(2) \to R$ is of the form 
$(d_{1},d_{2})\mapsto a+ d_{1}\cdot b_{1}+d_{2}\cdot b_{2}$
for unique $a$, $b_{1}$ and $b_{2}$ in $R$.}

\medskip 

 \noindent One may have deduced  Axiom 2 from Corollary \ref{31x}, 
{\em provided} one knew that any function $D(2)\to R$ may be extended 
to a function $D\times D \to R$. But this is not automatic - rather, 
this is guaranteed by the Axiom 2.

Of course, there are similar axioms for $n=3,4, \ldots$, using $$D(n):= 
\{ (d_{1}, \ldots ,d_{n})\in R^{n}\mid d_{i}\cdot d_{j}=0 \mbox{ for 
all } i,j=1, \ldots n \}.$$  
In short form,  a general Axiom 2 says: {\em Any map $D(n) \to R$ extends uniquely to an affine map 
$R^{n} \to R$.}

\medskip

Another  instantiation of the axiom scheme gives the following 
Axiom (we shall not use here): Let $D_{2}:= \{ x\in R\mid x^{3}=0\}$.
Then {\em every function $f:D_{2}\to R$ is uniquely of the form $x\mapsto a_{0} + a_{1}\cdot x 
 + a_{2}\cdot x^{2}$}, or: every $f:D_{2}\to R$ extends uniquely to a 
polynomial function $R\to R$ of degree $\leq 2$.  This axiom corresponds to the 3-dimensional Weil 
algebra ${\mathbb Q}[X]/(X^{3})$. More generally, let 
$D_{k}(n):=\{(x_{1}, \ldots ,x_{n})\in R^{n}\mid \mbox{all products 
of $k+1$ of the $x_{i}$s is } 0\}$. Then {\em every function $D_{k}(n) \to 
R$ extends uniquely to a polynomial function $R^{n}\to R$ of degree 
$\leq k$}. -- The polynomial functions occurring here are the Taylor 
polynomials at $0\in R$)  (resp.\ at $(0,\ldots,0)\in R^{n}$) of $f$.

As a final example of an instantiation of the axiom scheme, let $D_{L}\subseteq 
R^{2}$ be given by
$$D_{L}:= \{(x_{1},x_{2})\in R^{2}\mid x_{1}^{2}=x_{2}^{2} \mbox{ and 
} x_{1}\cdot x_{2}=0\}.$$ Then the following axiom is likewise an 
instantiation of the axiom scheme: {\em every function $f:D_{L}\to R$ is of 
the form $f(x)=a+b_{1}\cdot x_{1} + b_{2}\cdot x_{2} + c\cdot 
(x_{1}^{2}+x_{2}^{2})$}. The $c$ occurring here can then be seen as 
(one fourth of) the Laplacian $\Delta (f)$ of $f$ at $(0,0)$.
Note that $D(2) \subseteq D_{L} \subseteq D_{2}(2)$. The space $D_{L}$ corresponds 
to a certain 4-dimensional Weil algebra; see also \cite{SGM} 8.3.

\subsection{Envelopes again}\label{EAx}
This section is to ``justify'' in classical terms the correctness of our description 
of envelopes in terms of characteristics, as in Section \ref{CEx}. For sim\-pli\-city, we con\-sider a 
1-parameter family of (un\-parame\-trized) curves $S_{t}$ in $R^{2}$.
We assume that there is some smooth function $F(x,y,t)$ such that the 
$t$th curve $S_{t}$ is given as the zero set of $F(-,-,t)$. We then prove 
that the classical analytic ``discriminant'' description of the 
characteristics and the envelope agrees with the synthetic/geometric  
one which we 
have given; but note that our description is coordinate free, so in 
particular, it follows that the constructed envelope is independent 
of the analytic representation.  
 To say that $(x,y)$ belongs to the $t_{0}$ characteristic is by 
 the synthetic definition to say that  $F(x,y,t_{0}+d)=0$ for all 
$d\in D$ (the neighbours $t$ of $t_{0}$ are of the form $t_{0}+d$). 
Equivalently, by Taylor expansion,
\begin{equation}\label{partialx}F(x,y,t_{0})+ d\cdot \partial 
F/\partial t (x,y,t_{0})=0,\end{equation}
for all $d\in D$. By the principle of cancellation of universally 
quantified $d$s, this is equivalent to the conjunction of the two 
equations \begin{equation}\label{discrx}F(x,y,t_{0})=0 \mbox { \quad and  \quad }\partial 
F/\partial t (x,y,t_{0})=0,\end{equation} which is how the $t_{0}$ 
characteristic, and hence the envelope,  may be  
described by the discriminant method. 

However, Courant gives an example (\cite{Courant}, Example 10 in III.3) to show that ``the envelope need 
not be the locus of the points of intersection of neighbouring\footnote{The word 
``neighbouring'' here is not in the sense of the $\sim$ neighbour 
relation that we are using, in fact, it rather means: distinct.} 
 curves'', in other words, the ``{\em non-rigourous but intuitive}'' 
 description of characteristics suggested in loc.cit.,  is not only 
 non-rigourous, it is furthermore 
{\em wrong}. (So implicitly: don't believe in geometry!). The example is 
the following. Consider the family of curves in the plane given by 
$F(x,y,t)= y-(x-t)^{3}$. (This is the curve $y=x^{3}$, together with 
all its horizontal translates.) We leave to the reader to prove that the 
characteristic set at parameter value $t_{0}$ (as calculated by 
(\ref{discrx})) is the subset 
$\{(t_{0}+D , 0)\}$ of the $x$-axis; so the envelope is the $x$-axis.
Whereas the ``limit intersection point'' idea does not work here, since 
(to quote Courant) ``no two of these curves intersect each other''. 

\subsection{Defining $\sim$ in terms of $R$?}
We have already postulated that $x\sim y$ in $R$ means $y-x \in D$, 
or $(y-x)^{2}=0$.
A (first order) neighbour relation $\sim$ on any object $M \in R$ can  
be defined by
\begin{equation}\label{hier1x}x\sim y \mbox{ in $M$ }\mbox{ iff } \alpha (x) \sim \alpha (y) \in 
D \mbox{ for all } 
\alpha: M\to R.\end{equation}
So $\sim$ is, for all objects $M$,  defined in terms of the scalar 
valued functions on $M$. Trivially,  any map $M' \to M$ preserves 
$\sim$. This is the ``contravariant'' or ``weak'' way of defining 
$\sim$. There is also a ``covariant'' or ``strong'' way of defining 
it, see \cite{SGM} p.\ 31. For good spaces, like $R^{n}$, they 
coincide. (The weak determination is not adequate in algebraic 
geometry, since projective space, and other important geometric 
objects,  only admit constant scalar valued functions. So one must 
here replace the consideration of scalar-valued functions by {\em locally 
defined} scalar valued functions, and for this, one needs some notion 
of ``local'', as alluded to in Section \ref{AiGx}.)

For the weak determination of $\sim$ , we can  identify the monad $\M (\underline{0})$ around the origin  in 
$R^{n}$:
\begin{prop}\label{32x} We have $\M (\underline{0}) = D(n)$. 
\end{prop}
(For $n=1$, this was postulated.) Let us prove it for $n=2$. We have already seen in (\ref{protxx}) that
$\M (\underline{0}) \subseteq D(2)$. For the converse, we have to 
consider an arbitrary map $\alpha : R^{2}\to R$ and prove that 
$(d_{1}, d_{2}) \in D(2)$ implies $\alpha (d_{1}, d_{2})^{2}=0$. By Axiom 2,
$\alpha (d_{1}, d_{2})= a + b_{1}\cdot d_{1} + b_{2}\cdot d_{2}$, and 
so $\alpha (d_{1}, d_{2})- \alpha (0,0)= b_{1}\cdot d_{1} + 
b_{2}\cdot d_{2}$, which has square 0 since $d_{1}^{2}=d_{2}^{2}= 
d_{1}\cdot d_{2}=0$.

\medskip

 \noindent {\bf Remark.} Since $D(2)$ is strictly smaller than $D\times 
D$, we therefore also have  that $(d_{1}, d_{2})\sim (0,0)$ is 
stronger than the conjunction  of $d_{1}\sim 0$ and $d_{2}\sim 0$.

\medskip

 The  Principle in the beginning of Chapter  
\ref{NCEGx} may now be proved algebrai\-cally:  we may assume that 
coordinates are chosen so that the considered common point $x\in 
l_{1}\cap l_{2}$ is $(0,0)$, and that $l_{1}$ and $l_{2}$ are graphs 
of the functions $x\mapsto b_{1}\cdot x$ and $x\mapsto b_{2}\cdot x$.
We must prove that $b_{1}=b_{2}$. For $d\in D$, we have 
$(d,b_{i}\cdot d) \in \M(x) = D(2)$. So by assumption, for all $d\in 
D$, we have 
$(d,b_{1}\cdot d) \in \M(x)\cap l_{1}= \M (x)\cap l_{2}\subseteq 
l_{2}$, so for all  $d\in D$, we have
$b_{1}\cdot d = b_{2}\cdot d$; cancelling the universally quantified 
$d$ then gives $b_{1}=b_{2}$.

\subsection{Contravariant and covariant hierarchy}
The polynomial function $x: R \to R$ vanishes at $0$; one also says 
that it vanishes to {\em first} order at $0$, and that $x^{2}$ 
vanishes to {\em second } order at $0$, etc.; more generally, $f:R 
\to R$ vanishes to second order at $0$ if it may be written 
$f(x)=x^{2}\cdot g(x)$ for some function $g:R\to R$.  Similarly for 
$k$th order vanishing. It generalizes to ``order of vanishing''  of $f:M 
\to R$ at a 
point $a\in M$. Note that $k$th order vanishing 
is a {\em weaker} condition 
than $(k+1)$st order vanishing.  This (essentially classical) 
hierarchy of scalar valued functions (quantities) is to be 
compared with the hierarchy of neighbours, applicable to points of 
spaces $M$, where $k$th order neighbour is a {\em stronger} condition than 
$(k+1)$th order neighbour. The neighbour relations are covariant notions, 
applicable to {\em points (ele\-ments)} of spaces (the assumed neighbour relations $\sim_{1}$,$\sim_{2}$,\ldots , are 
preserved by mappings, and thus are {\em covariant}); 
the order-of-vanishing is a contravariant notion, applicable to 
{\em quantities} on $M$ , i.e.\ to $R$-valued functions $M\to R$.

The notions are related as follows, for $a$ and $b$ in $M$:  $a\sim _{k}b$ iff 
for any quantity $f:M\to R$ vanishing to $k$th order at $a$, we have 
$f(b)=0$;  this is, for $k=1$, just a reformulation of (\ref{hier1x}). Recall 
that 
$a\sim _{k}b$ on $R$ is defined in terms of $(a-b)^{k+1}=0$, i.e.\ in 
terms of order of {\em nilpotency}.  

A 
classical formulation, in certain contexts, is that we can ``ignore'' quantities of higher
 order, in comparing $a$ and $b$: ``{\em Dabei sehen wir von unendlich kleinen 
Gr\"{o}ssen h\"{o}here Ordnung ab.}'' (``Here, we ignore infinitely small 
quantities of higher order.''), \cite{BT} p.\ 523. In rigourous 
mathematics, one cannot  
``ignore'' anything except $0$. But one can certainly consider nilpotent elements in 
rings. Thus,
an explicit theory of infinitesimals came in through the back door, 
namely from algebraic geometry:

\subsection{Wisdom from algebraic geometry}\label{wisdx}
The development leading to the modern formulations of SDG began in 
French algebraic geometry in the mid 20th Century by Grothendieck and 
his collaborators, with the notion (and category!) of 
{\em schemes}, as a generalization of the notion of algebraic 
varieties 
(over a field $k$, say). 

In particular, the category $\E _{k}$ of affine schemes over $k$ is 
the by definition the dual of the category $\A_{k}$ of commutative 
$k$-algebras, suitably size-restricted, say: of finite presentation. 
The algebras are allowed to have nilpotent elements. 
Such algebra $A$ is  seen as the ring of scalar valued functions 
on the scheme $M$ (geometric object, ``space'') which it defines. One 
writes $M=\Spec (A)$. 
(The ``scalars'' $R$ is the scheme represented by $k[X]$.)
Then the 
algebra
$A\otimes A$ defines the space $M\times M$. If $I$ is the kernel if 
the multiplication map $A\otimes A \to A$, then $(A\otimes A)/I\cong A$. 
Consider the ideal $I^{2}\subseteq I$. The $k$-algebra $(A\otimes 
A)/I^{2}$ gives $M_{(1)}$, the first neighbourhood of the diagonal of 
$M$. So $$M_{(1)}:= \Spec ((A\otimes A)/I^{2}).$$
The quotient map $(A\otimes 
A)/I^{2} \to (A\otimes A)/I \cong A$ defines, in the category of 
schemes, the diagonal $M \to M_{(1)}$.

Note that $I/I^{2}\subseteq (A\otimes A)/I^{2}$ consists of elements 
of square 0. It is in fact the module of {\em K\"{a}hler 
differentials of $A$};  $(A\otimes A)/I^{2}$ is the ring of scalar 
valued functions on $M_{(1)}$, and the submodule $I/I^{2}$ consists 
of  those 
functions that vanish on the diagonal $M \subseteq M_{1}$, i.e.\ the 
combinatorial scalar valued 1-forms, in the sense of Section 
\ref{DFx}. (K\"{a}hler introduced these differentials already in the 
1930s.)   

\medskip

The simplest  scheme which is not a variety 
is $D$, the affine scheme given by $k[\epsilon ]$, the ring of dual 
numbers
 over $k$. The underlying variety of $D$ has just one global point, 
since $k[\epsilon ]$ has only one prime ideal, namely $(\epsilon )$. 
Geometrically, $D$ is a ``thickened'' version of its unique global point.
Mumford (\cite{Mumford} p.\ 338) describes $D$  as ``a sort of 
disembodied tangent vector'', meaning that a map $D \to X$ may be 
identified with a tangent vector to $X$, for any scheme $X$.

The relationship between the infinitesimal objects like $D$, and the 
neighbourhoods of diagonals may be exemplified by the isomorphism
$$R_{(1)}\cong R\times D,$$
given by $(x,y) \mapsto (x, y-x)$, for $x\sim y$ in $R$.

The crucial step in the formation of contemporary SDG was when 
Lawvere in 1967 combined this consideration of a ``tangent vector 
representor'' $D$  with the idea of Cartesian closed category $\E$. Thus, for 
any object $X$ in $\E$, $X^{D}$ is then the {\em object} (space) in $\E$ of all 
tangent vectors to $X$, in other words, it is the (total space of the) 
tangent bundle $T(X)\to X$.

To put this relationship into axiomatic form is most conveniently done 
by assuming a ring object $R$, and describing $D$ in terms of $R$, (as 
is done in Section \ref{CGx}). There is a more radical approach, 
advocated by Lawvere in \cite{L:D}, where the ring $R$ is to be 
constructed out of an infinitesimal object $T$ (``an instant of 
time'') (ultimately then 
proved to be isomorphic to $D$); see also \cite{CC} 5.3.

\section{Models of the axiomatics}\label{MAx}
For an axiomatic theory, models are useful, but not crucial. 
Euclidean geometry has been useful for more than two thousand years. 
When exactly was a model for it presented? Did it have to wait for the real 
numbers, or at least some subfields of it?
Models are useful, - they may guide the intuition, and prevent inner 
contradictions. This also applies to SDG. The models for SDG come in two main 
groups: arising from algebraic geometry, and from classical differential geometry over 
${\mathbb R}$, respectively (and in fact, SDG serves to make explicit 
what the two groups have in 
common).

Models for the axiomatics of Section \ref{NMCx} may be 
built on the basis of some of the models of SDG mentioned above; see 
\cite{MSSDG}.

\subsection{Algebraic models}
The category $\E _{k}$ of affine schemes over a commutative ring $k$ (i.e., the 
dual of the category of (finitely presented, say) commutative $k$-algebras) 
is a model\footnote{If 
$2$ is not invertible in $k$, there are things that work differently.}, with $k[X]$ as $R$. 
$\E _{k}$ is not quite Cartesian closed, but at least the scheme 
corresponding to 
$k[\epsilon ]$ (or to any other Weil algebra) is exponentiable. 
The set valued presheaves $\hat{\E _{k}}$ on $\E _{k}$ is a full 
fledged topos model (with $R$ 
represented by $k[X]$). The topos $\hat{\E _{k}}$ is of course the 
same as the category of covariant functors from the category of 
(finitely presntable) commutative 
$k$-algebras to sets, and $R$ is in this set up just the forgetful 
functor, since $k[X]$ is the free $k$-algebra in one generator.

Many of the subtoposes of  $\hat{\E _{k}}$ are likewise models; 
passing to suitable subtoposes, one may force $R$ to have further properties; one 
may for instance force $R$ to become a local ring; the subtopos 
forcing this is also known as the Zariski Topos. These toposes are 
explicitly the main categories studied in \cite{DG}. 

\subsection{Analytic models based on ${\mathbb R}$}
There is of course a special interest in models $(\E ,R)$ which contain the 
category $Mf$ of smooth manifolds as a full subcategory, in a way which 
preserve  known constructions and concepts from classical 
differential geometry. So one wants a full and faithful functor 
$i: Mf \to \E$, with $i({\mathbb R})=R$. Also {\em transversal} 
pull-backs should be preserved, and $i(T(M))$ should be $i(M)^{D}$. 
The properties of such a functor $i$ has been axiomatized by Dubuc \cite{Dubuc}) 
under the name of ``well adapted model for SDG''; see also 
\cite{PWAM}. The book \cite{MR} is 
mainly devoted to the construction and study of such models.

The earliest well-adapted model  (constructed 
by Dubuc \cite{Dubuc})  is  one now known as the 
``Cahiers topos''. It can be proved to contain the category of 
convenient vector spaces (with smooth maps between them) as a full 
subcategory, in a way which preserves the Cartesian closed structure, 
cf.\ \cite{KR1}, \cite{KR2}. A more advanced topos $\GG$, now called 
the``Dubuc topos'',  
 \cite{Dubuc2}, even supports some ``Synthetic Differential Topology'', cf.\ 
\cite{SDT}.

A  main tool in the construction of analytic models is to take heed 
of the wisdom of algebraic geometry, but replacing the algebraic 
theory (in the sense of Lawvere) ${\mathbb T}$ of 
commutative rings with the richer algebraic theory ${\mathbb T}_{\infty}$, whose 
$n$-ary operations are not only the real polynomial functions, but 
{\em all} the smooth maps ${\mathbb R}^{n} \to {\mathbb 
R}$. It contains the theory of commutative rings as a subtheory, 
since a polynomial in $n$ variables defines a smooth function in $n$ 
variables. The  theory ${\mathbb T}_{\infty}$, and its importance for the project of 
categorical dynamics, was already in Lawvere's seminal 1967 lectures.

Note that any smooth manifold $M$ gives rise to an algebra for this 
theory, namely $C^{\infty}(M)$, the ring of smooth ${\mathbb 
R}$-valued functions on $M$. We may think $M$ as a ``reduced'' affine 
scheme corresponding to the ring $C^{\infty}(M)$, and then mimick the 
construction (described above) of set valued presheaves on affine 
schemes, and subtoposes thereof.  But note that also  ${\mathbb R}[\epsilon ]$ (and 
all other Weil algebras over ${\mathbb R}$) are algebras for 
${\mathbb T}_{\infty}$, and define (non-reduced) affine schemes.

Modules of K\"{a}hler differentials for  algebras for ${\mathbb 
T}_{\infty}$ were studied in \cite{DK}  

(If one takes just the category of smooth manifolds (with open 
coverings) as site of 
definition for a topos, one gets a topos already considered in SGA4, 
under the name of ``the smooth topos''; it contains the category of 
diffeological spaces as a full subcategory, but lacks the 
infinitesimal objects like $D$. These categories are models for the 
Fermat-Reyes axiomatics. See \cite{SDG} Exercise III.8.1).

\section{New spaces} 
Except for the ``infinitesimal'' spaces like 
$D_{k}(n)$, the present account does not do justice to 
the {\em new} spaces which have emerged through the development of 
SDG. In particular, it has not capitalized on the unproblematic way 
in which function spaces exist in this context, by Cartesian 
closedness of $\E$. These function spaces 
opens the 
door to a synthetic treatment of calculus of variations, continuum 
mechanics, infinite dimensional Lie groups, \ldots . For such spaces, 
the neighbour relation (which has been my main focus here) is  more 
problematic, however,  and is not 
well exploited. Instead, one uses the (classical) method of encoding the 
infinitesimal information of a space $X$ in terms of its tangent 
bundle $T(X)= X^{D}$, rather than in terms of $X_{(1)}$ (first 
neighbourhood of the diagonal). Notably Nishimura has pushed the 
SDG-based theory far in this direction, cf.\ e.g.\ \cite{N}.

Another type of new spaces come from the observation that the functor 
$(-)^{D}$ in many of the models has a right adjoint, $(-)^{1/D}$ 
(Lawvere's notation, ``fractional exponent''); the spaces $M^{1/D}$ 
are re\-mi\-ni\-scent of Eilenberg-Mac Lane spaces. There is some 
discussion of them in \cite{SDG} I.20, in \cite{KR3}, and in 
\cite{L98}. 

Finally, the notion of {\em jets}, and the jet bundles, as considered 
by Ehresmann in the 1950s, form, on the one hand, one of the sources 
for SDG as presented here; on the other hand, the SDG method  makes 
the consideration of jets and jet bundles simpler, since SDG makes 
the notion of jet {\em representable}, in the sense that a $k$-jet at 
$x\in M$, with values in $N$, is here simply a map $\M _{k}(x)\to N$, 
rather than an equivalence class of maps $U\to N$ (where $x\in U$). (In 
\cite{SDT}, the notion of {\em germ} of a map is likewise 
representable.)

When jets are representable, Ehresmann's 
theory of differentiable group\-oids, as carrier of a general theory of 
connections,  admits some simpler formulations:

\subsection{Connections in fibre bundles and groupoids}\label{Connx}

For the present purpose, a {\em fibre bundle} over a space $M$ is 
just a map $\pi: E\to M$. (When it comes to proving things, one will 
need good exactness properties of $\pi$, like being an effective 
descent map, or being locally a  projection $F\times M \to M$.) Then a combinatorial connection 
in the bundle $E\to M$ is 
an action $\nabla$ of $M_{(1)}$ on $E$, in the sense that  $(x,y)\in 
M_{(1)}$  and $e\in E_{x}$ define an element 
$\nabla(x,y)(e)$ in $E_{y}$. (Here, $E_{x}:= \pi^{-1}(x)$, and 
similarly for $E_{y}$.) One requires the normalization condition 
$\nabla(x,x)(e)=e$. For good spaces, it then follows that 
$\nabla(y,x)\nabla (x,y)(e) =e$. The notion of affine connection $\lambda$  
considered above is a special case: the bundle $E\to M$ is in this 
case the first projection $M_{(1)} \to M$, and $\nabla (x,y)(z)= 
\lambda (x,y,z)$. If $E \to M$ is a vector bundle, say, a {\em linear} 
connection is a connection $\nabla$  where the map $\nabla (x,y)(-): 
E_{x}\to E_{y}$ is linear for all $x\sim y$. For good spaces $M$, 
linear connections in the tangent bundle $T(M)\to M$ contain exactly 
the same information as affine connections $\lambda$ on $M$.

There is also a  
notion of connection $\nabla$ in a {\em groupoid}  $\Phi 
\rightrightarrows M$. (This is 
closely related to the notion of {\em principal connection} in a 
principal fibre bundle $P\to M$; in fact, such $P$ defines, according 
to C.\  
Ehresmann, a groupoid $PP^{-1}\rightrightarrows M$, and a principal 
connection in $P\to M$ is then the same data as a groupoid connection in 
$PP^{-1}\rightrightarrows M$.) Recall 
that a groupoid $\Phi \rightrightarrows M$ carries a reflexive 
symmetric structure: the reflexive structure picks out for every 
$x\in M$ the identity arrow at $x$, and the symmetric structure 
associates to an arrow $f:x\to y$ its inverse $f^{-1}:y \to x$. Then 
a connection in $\Phi \rightrightarrows M$ is simply a map 
$M_{(1)}\to \Phi$ preserving (the two projections to $M$ and) the reflexive and symmetric structure,
$$\nabla (x,x)= \id_{x}\mbox{\quad  and \quad  } \nabla(y,x)= \nabla (x,y)^{-1}$$
for all $x\sim y$. 

Given a bundle $E\to M$ in $\E$. If $\E$ is locally Cartesian closed, 
one may form the groupoid $\Phi \rightrightarrows M$ where the arrows 
$x\to y$ are the invertible maps $f: E_{x}\to E_{y}$. Then a connection on 
$E\to M$, in the 
bundle sense, is equivalent to a connection, in the groupoid sense, 
of this groupoid $\Phi \rightrightarrows M$. If $E\to M$ is a vector 
bundle, there is a subgroupoid of $\Phi \rightrightarrows M$ 
consisting of the {\em linear} isomorphisms $E_{x}\to E_{y}$ (this groupoid 
deserves the name $GL(E)$). Similarly if $E\to M$ is a group bundle, 
or has some other fibrewise structure.

The groupoid formulation of the notion of connection is well suited 
to formulate algebraic properties, like curvature. We may observe 
that the curvature, as described in Section \ref{ACx} for affine 
connections $\lambda$, is purely groupoid theoretical. Thus if $x, y, 
z$ form an infinitesimal 2-simplex in $M$, it makes sense to ask 
whether $\nabla (x,y)$ followed by $\nabla (y,z)$ equals $\nabla 
(x,z)$, or better: consider the arrow $R(x,y,z): x\to x$ given as the composite 
(composing from left to right)
$$R(x,y,z):=\nabla (x,y).\nabla (y,z).\nabla 
(z,x) \in \Phi (x,x).$$
This is the {\em curvature} of $\nabla$, more precisely, the 
curvature $R$ is a 
combinatorial 2-form with values in the group bundle $gauge (\Phi)$ of vertex groups 
$\Phi (x,x)$ of $\Phi$. Now the connection $\nabla$ in $\Phi$ gives rise to a 
connection $ad\nabla$ in the group bundle $gauge (\Phi)$: $ad \nabla 
(x,y)$ is the (group-) isomorphism $\Phi (x,x) \to \Phi (y,y)$ 
consisting in conjugation by $\nabla (x,y): x\to y$. This conjugation we  
write $(-)^{\nabla(x,y)}$. In terms of this, we have an identity, which 
deserves the name the {\em Bianchi} identity for (the curvature $R$ 
of) the connection 
$\nabla$; namely for any infinitesimal 3-simplex $(x,y,z,u)$, we have
\begin{equation}\label{Bianx}id_{x}= R(yzu)^{\nabla (y,x)}. 
R(xyu).R(xuz).R(xzy),\end{equation}
verbally, {\em the covariant derivative of the $gauge (\Phi)$ valued 2-form  $R$, 
with respect to the connection $ad \nabla$ in the group bundle, is 
``zero''}, i.e.\ takes only identity arrows as values.

The proof of (\ref{Bianx}) is trivial, in the sense that it is a case 
of Ph.\ Halls 14-letter identity, which holds for any six elements in 
a group, or for the six arrows of a tetrahedron-shaped diagram in a 
groupoid; here, the six arrows are the $\nabla(x,y)$, $\nabla (x,z)$, 
\ldots , $\nabla (z,u)$ in $\Phi$. See \cite{SGM}, and see  \cite{CCBI} 
for how this implies the classical Bianchi identity for linear 
connections in vector bundles.

\section{The role of analysis}

\subsection{Analysis in geometry?}\label{AiGx}The phrase {\em analytic} geometry may 
be used in the wide sense: using coordinates and calculations. In 
this sense,  SDG as presented here quickly becomes analytic (e.g.\ 
the basic axiomatics is formulated in such terms, as expounded in 
Chapter \ref{CGx}). But the more common use of the phrase ``analytic'' 
is that {\em limit} processes  and topology are utilized.

Ultimately, topology and limits in real analysis have their origin in 
the strict order relation $<$ on ${\mathbb R}$. Then the partial 
order $\leq$ is defined by $x\leq y$ iff  $\neg (y<x)$. The elements
 in ${\mathbb R}_{>0}$ are invertible. In SDG, it is also 
natural to have  an order $<$ on $R$, given primitively, or in terms of the 
algebraic structure of $R$. (In well adapted models $i: Mf \to \E$, the relation $<$ 
is definable in terms of the inclusion of the smooth manifold 
${\mathbb R}_{>0}$ into $ {\mathbb R}$, which by the embedding $i$ 
defines a subobject $R_{>0}\subset R$, out of which  a strict 
order $<$ can be defined.) Nilpotent elements $d$ in a 
non-trivial ring
cannot be invertible. It follows, for any nilpotent $d$, that $d\leq 0$, 
and hence also $-d\leq 
0$ (since also $-d$ is nilpotent). So $0\leq d \leq 0$. So if $\leq$ 
were a partial order (not just a 
preorder), this 
would imply that any nilpotent $d$ is $0$, which is incompatible with SDG. Thus, in SDG, $\leq $ 
is only a preorder, not a partial order. For preordered sets, a supremum 
is not uniquely defined; to have a {\em unique} number as 
supremum, one needs a partial order.

This is one reason why limit processes are not used in SDG, at the 
present stage. 

Topology  comes in 
play e.g.\ when formulating  statements about {\em local} existence 
of, say, solutions to particular  differential equations. `Local' 
refers to some topology on a given object, and in SDG, there may be
several natural choices, cf.\ in particular the recent \cite{SDT}. The finest topology on an object (space) $X$ 
is, in the context of SDG, the one where the open subsets are those 
$U\subseteq X$ which are closed under the neighbour relation $\sim$. 
For instance,  a local solution $f$ (for this fine topology) for a differential equation 
$f'=F(x,y)$ amounts to a formal power series solution, and is 
therefore cheap. More serious existence statements are when stronger 
topologies, like the ``intrinsic Zariski topology'', are involved : an 
subset $U\subseteq X$ is open if it is of the form $f^{-1}(R^{*})$, where 
$R^{*}\subseteq R$ consists of the invertible elements; or if it is 
of the form $f^{-1}(P)$, where $P \subseteq R$ 
consists in the strictly 
``positive'' numbers - which then in turn have to be described or 
assumed; see \cite{Mcl}, for some results in this direction.

SDG does not {\em prove} basic integration results, and even 
the formulation of such results does not come for free. Advances in 
this direction exists, in what is now called Synthetic Differential 
Topology. It builds on SDG, and its main model is 
the Dubuc topos $\GG$; see  \cite{MR} Chapter III, and 
notably \cite{SDT}, where also a synthetic theory of singuarity 
theory is considered.

The most basic integration result is the (essentially unique) 
existence of anti-derivatives: for 
$f:R\to R$, there exists $F:R\to R$ with $F' =f$. In an axiomatic 
development, this has to be 
taken as an {\em axiom}, -- one that actually can be proved to hold in all the 
significant topos models $(\E, R)$ for SDG. Similarly for many other 
basic results, like a suitable version of the intermediate value 
theorem.\footnote{Significantly, the version valid in significant 
SDG models applies to functions $f$ with a {\em transversality} condition,  like 
$f' >0$, -- like in constructive analysis.}  Thus, full fledged 
analysis in axiomatic terms, incorporating SDG, quickly becomes 
overloaded with axioms, and is better developed as a {\em 
descriptive} theory, describing what actually {\em holds} in {\em 
specific} models $(\E, R)$. This is the approach of \cite{MR} which 
significantly has the title ``Models for Smooth Infinitesimal 
Analysis'' (although also a full-fledged axiomatic theory is 
presented in loc.cit., Chapter VII). Note 
that the term ``smooth'', in so far as SDG is concerned, is a void term, 
since unlimited differentiability is automatic in this context; and 
``smooth implies continuous'' (equivalently, ``all maps are 
continuous'') is a {\em Theorem} in the good 
well-adapted models, see e.g.\ Theorem III.3.5 in \cite{MR}. 

I prefer not to think of SDG as a monolithic global theory, but as a {\em 
method}  to be used locally, in situations  where 
it provides insight and
simplification of a notion, of  a construction, or of an argument. 
The assumptions, or axioms that are  needed, may be taken from the valuable 
treasure chest of real anaysis.

Thus, the very construction of well adapted models $Mf \to \E$ 
depends on the theory ${\mathbb T}_{\infty}$ whose $n$-ary operations are the smooth 
functions ${\mathbb R}^{n}\to {\mathbb R}$, so that e.g.\ the exponential function
 $\exp: {\mathbb R}\to {\mathbb R}$, or the trigonometric functions, 
are ``imported'' from the treasure 
chest (here, imported from Euler, say, much prior to the rigourous formulation of 
limit processes). In the context of SDG, it is possible to  introduce 
existence of, say, these particular transcendental functions 
axiomatically, by  functional equations, or 
by  differential equations.  This is what the Calculus 
Books in essence do.

\subsection{Non-standard analysis ?}
Non standard analysis (NSA) is another theory where the notion of 
infinitesimals has an explicit and well defined status. Therefore, 
 one sometimes asks whether there is some relationship between SDG 
and NSA. There is 
very little relationship; NSA is a descriptive, not an axiomatic, theory, 
dealing (at least in so far 
as differential geometry goes) with the real 
number field ${\mathbb R}$, and crucially capitalizing on  its Cauchy  
completeness, since it is crucial that {\em  every (bounded) non-standard real number $\in {\mathbb 
R}^{*}$ has a unique standard part}. This is another expression of  
the completeness of the real number system. In this sense, NSA is a 
reformulation, with a richer vocabulary, of standard real analysis, and can, as 
such, cope with things defined in terms of limits, like definite 
integrals in terms of Riemann sums, say; SDG cannot do this, at best, 
it can introduce some integration by axioms, cf.\ the remark on the 
Frobenius integration Theorem in Section \ref{GDx}.

In NSA, one has a neighbour notion for elements in ${\mathbb R}^{*}$; 
it is an {\em equivalence} relation, and the equivalence classes are called 
{\em monads} -- a term which SDG has imported; but in SDG it is 
crucial that the neighbour relations are not transitive, and come in a hierarchy: first order, 
second order, \ldots, (hence first order, second  order, \ldots , 
$k$th order monads $\M _{k}(x)$, \ldots ), and this comes closer to important
aspects of mathematical practice, where notably the first order 
neighbour relation takes most of the work on its shoulders, and has 
been the sole concern in this note. (The second order  monads in SDG play a 
role when discussing e.g.\ dynamic or metric notions; thus a (pseudo-) Riemannian metric 
may be defined in terms of a $R$-valued functions $f(x,y)$ defined for 
$x\sim_{2}y$, and with $f(x,y)=0$ if $x\sim_{1}y$; see \cite{SGM}.)

NSA can also be axiomatized, but this amounts essentially to 
axiomatizing a further structure (an endo-functor) on the {\em category} of (discrete) sets 
\cite{KM}, or a further primitive predicate in axiomatic (Zermelo 
Fraenkel) set 
theory \cite{Nelson}.

\section{The continuum and the discrete}\label{Phil}
An historically important  problem in (the philosophy of) mathematics 
is the problem of understanding the nature of the continuum, and its 
relationship to the discrete.
 Is the continuum just a discrete set of points? (and motion 
therefore impossible, according to Parmenides). In contrast,  in Euclidean 
geometry,  {\em line} (line segment) was a primitive notion, and was not just 
the set of points in it. (And {\em time} was not a set of instants.) 
Even a contemporary geometer like Coxeter 
makes the distinction between a line and the ``range of points'' on 
it, cf.\ \cite{Cox} p.\ 20. 

The principal side of the 
contradiction between continuum and discrete was, historically, the 
continuum. With the full arithmetization of the continuum, in the 
hands of, say, Dedekind, with the construction of the real number 
system ${\mathbb R}$, the continuum was reduced to a set of points, 
and the cohesion of the continuum was reduced to a topology on this 
point set. For mainstream differential geometry 
synthetic axiomatic considerations 
became, in principle, redundant. Everything became reducible to real 
analysis.

Synthetic differential geometry refuses to take this  one-sided 
reductionist view. (For one reason. ${\mathbb R}$, as a point set 
(set of global points), has 
no non-trivial nilpotent elements $d$.)
Rather, SDG learns from (and possibly contributes to) {\em analyzing} the relationship between the 
continuum  and the discrete. Such analysis typically has the form of 
a functor $\gamma_{* }: \E \to \Ss$, with $\Ss$ some category of discrete 
sets, and with $\E$ some category of spaces with some kind of 
cohesion\footnote{a situation axiomatized  in Lawvere's \cite{LC}
 ``Mengen'' vs.\ ``Kardinalen'', and further elaborated in papers by 
 Lawvere and by Menni, cf.\ e.g.\ \cite{Menni} and \cite{LaMe}.}
Preferably, both $\E$ and $\Ss$ are toposes, and $\gamma_{*}$ a geometric 
morphism, associating to a space $X \in \E$ its set of (global) 
points. The left adjoint $\gamma^{*}$ of $\gamma*_{}$ is a full 
embedding, so that discrete spaces form a full subcategory of $\E$. 
An example of such $\E$-$\Ss$-pair is with $\E$ the topos of 
simplicial sets, with $\gamma_{*}(X)$ the set of $0$-simplices  (= 
global points) of $X$. 
This example is relevant to algebraic topology (cf.\ e.g.\ \cite{May}), not to differential 
geometry, but it illustrates a phenomenon which is crucial also for SDG: 
namely that there are non-trivial objects with only one global point 
(e.g.\ in the topos of simplicial sets: the simplicial $n$-sphere 
$\Delta (n)/\dot{\Delta}(n)$) 
- just like $D$ in SDG has $0$ as the only global point.

A  well-adapted model $\E$ of SDG contains not only the the category 
of discrete 
manifolds (sets) as a full subcategory, but even the category of all smooth 
manifolds, in particular ${\mathbb R}$. By the fullness, ${\mathbb 
R}$, when seen in $\E$ (and there denoted $R$)  
does not acquire any new global points (unlike 
the ${\mathbb R}^{*}$ of NSA). But it does acquire new subobjects, - 
e.g.\ $D \subseteq R$. When we talk about general elements $d\in D$, we are 
therefore not talking about {\em global} points $1\to D \subseteq R$. 

A space is an object in a category of spaces (Grothendieck, Lawvere). 
So what ``is'' the space  ${\mathbb R}$? It depends on the category 
in which it is considered. In SDG, one considers ${\mathbb R}$ in  
certain (``well adapted'') toposes $\E$; ${\mathbb R}=R$ does not 
change, it is the ambient category which changes.

\section{Looking back}\label{Lbx}
The discovery, by Huygens in the 17th Century, of the notion of 
envelopes and their relatives, (leading to a theory of waves, 
isochrones, \ldots ), was coined in geometric terms, without essential reference (so far I know) to analytic considerations. When differential 
calculus, as we know it today, was developed, analytic methods became 
more dominant. A main treatise like Monge's in 1795 was entitled 
``{\em L'application de l'analyse \`{a} la g\'{e}om\'{e}trie}''. But 
this treatise of Monge's goes  also in the 
other direction: it forcefully uses geometric and synthetic reasoning 
for explaining the analytic  theory of first order PDEs of Lagrange,
 - a thread taken up later by 
Sophus Lie; this comprises in particular  the theory of {\em characteristics} of such PDEs, 
the curves, out of which the solutions of the PDE can be built. (They 
are built up from characteristics in the sense of Section 
\ref{CEx}, namely intersections of families of surface elements.) 
Lie's 1896 book \cite{BT} on contact geometry
has  a chapter called
{\em Die Theorie der partiellen Differentialgleichungen als Teil der 
Theorie der Fl\"{a}chenelemente}.
(Fl\"{a}chenelement = surface element = contact element $\M(b)\cap B$, as in Section 
\ref{NMCx}, or the sets $\M_{\approx}(x)$ of suitable codimension 1 
distribution, as in Section \ref{GDx}.)

In one of Lie's 
early articles on the theory of differential equations, he wrote:

\medskip 

{\em ``The reason why I have postponed for so long these investigations,
which are basic to my other work in this field, is essentially the following.
I found these theories originally by synthetic conside\-rations. 
But I soon realized that, as expedient [zweckm\"{a}ssig] the synthetic 
method is for discovery, as difficult it is to give a clear exposition 
on synthetic investigations, which deal with objects that till now have 
almost exclusively been considered analytically. After long vacillations, 
I have decided to use a half synthetic, half analytic form. 
I hope my work will serve to bring justification to the synthetic method 
besides the analytical one."}

\medskip

(From Lie's ``Allgemeine Theorie der partiellen Differentialgleichungen
erster Ordnung'', Math.\ Ann.\ 9 (1876); my translation.)

In spite of Lie's call for a synthetic language and logic, the 
differential geometry in the 20th 
Century became more and more analytic, and removed from the geometric 
intuition - at the time of Einstein, the ``d\'{e}bauche of indices'', and 
rules for how the coordinates transform, later on more abstract and 
coordinate free, but still somewhat un-geometric - as it must be when explicit 
infinitesimals (neighbour points) have to be avoided.

The editors of the present volume asked me to address the 
question about the ``advantages of SDG over other approaches \ldots ".
First of all, the neighbour notion, and  synthetic reasoning and 
concept formation 
 with 
it, is not an invention of present day SDG; it has been, and is, used again and 
again by engineers, physicists, by Sophus Lie (cf.\ the above 
quotation), by David Hilbert \cite{HCV}, and  (at least secretly) also by later 
mathematicians. However, explicit rules for such concept-formation, construction 
and reasoning
 have not been well formulated, and SDG is an attempt to provide such 
rules, so that the concepts, constructions and reasoning can be clearly 
communicated, and tested for rigour. What is the advantage of 
communication and rigour? It is not a question of ``advantage'', but 
a question of necessity.

\medskip

\noindent Anders Kock,\\
Dept.\ of Math.\ \\
University of Aarhus, Denmark\\
kock (at) math.au.dk



\small


 \begin{thebibliography}{99}
\bibitem{Be} W.\ Bertram, {\em Calcul diff\'{e}rentiel topologique 
\'{e}l\'{e}mantaire}, Calvage \& Mounet 2010.

\bibitem{BCS}R.\ Blute,  J.R.B.\ Cockett,
T.\ Porter, R.A.G.\ Seely, K\"{a}hler Categories, {\em Cahiers de 
Topologie et G\'{e}om\'{e}trie
Diff\'{e}rentielle} 52 (2011), 253-268.



\bibitem{BM} L.\ Breen and W.\  Messing, Combinatorial 
Differential Forms, {\em Advances in Math.\ } 164 (2001), 203-282.


\bibitem {SDT} M.\ Bunge, F.\ Gago, and A.M. San Luis, {\em Synthetic Differential Topology}, 
Cambridge University Press (to appear).

\bibitem{Busemann} H.\ Busemann, On spaces in Which Two Points 
Determine a Geodesic, {\em Trans.\ Amer.\ Math.\ Soc.\ } 54 (1943), 171-184.

\bibitem{CC}R.\ Cockett and G.\ Cruttwell, Differential Structure, 
Tangent Structure, and SDG, {\em  Appl.\  Categor.\ Struct.\ } 22 
(2014), 331-417.


\bibitem{Courant} R.\ Courant, {\em Differential and Integral 
Calculus} Vol.\ II, Blackie \& Son 1936.
\bibitem{Cox} H.S.M.\ Coxeter, {\em The Real Projective Plane}, 2nd 
ed., Cambridge University Press 1955.
\bibitem{DG} M.\ Demazure and P.\ Gabriel, {\em Groupes Alg\'{e}briques}, 
Tome I, Masson \& Cie/North Holland 1970.
\bibitem{Dubuc}  E.J.\ Dubuc, 
Sur 
les mod\`{e}les de la g\'{e}om\'{e}trie diff\'{e}rentielle 
synth\'{e}tique, {\em Cahiers de Top.\ et G\'{e}om.\ Diff.}\ 20 (1979), 
231-279.
\bibitem{Dubuc2}  E.J.\ Dubuc, $C^\infty$-schemes, {\em Amer.\ J.\ 
Math.\ } 103-104 (1981), 683-690.
\bibitem{DK}  E.J.\ Dubuc and A.\ Kock, On 1-form classifiers, {\em 
Communications in Algebra} 12 (1984), 1471-1531.

\bibitem{HCV}D.\ Hilbert and S.\ Cohn-Vossen, {\em Anschauliche Geometrie}, 
Grundlehren der mathematischen Wissenschaften 37, Springer 
Verlag 1932.



\bibitem{PWAM}A.\ Kock, Properties of well-adapted models for 
synthetic differential geometry, {\em Journ.\ Pure Appl.\ Alg.\ } 20 (1981), 
55-70.
\bibitem{SDG} A.\ Kock, {\em Synthetic Differential Geometry}, London Math. 
Soc.\ Lecture Notes Series 51 (1981); 2nd ed., London Math. 
Soc.\ Lecture Notes Series 333 (2006).

\bibitem{DFVG}A.\ Kock, Differential forms with values in groups, 
{\em Bull.\ 
Austral.\ Math.\ Soc.\ } 25 (1982), 357-386.


\bibitem{CTC}A.\ Kock, A combinatorial theory of connections, in 
{\em Mathematical Applications of Category Theory}, Proceedings 1983 
(ed.\ J.\ Gray),  A.M.S.\ Contemporary Mathematics Vol.\ 30 (1984).


\bibitem{KR1} A.\ Kock, Convenient vector spaces embed into the 
Cahiers topos, {\em Cahiers de Top.\ et G\'{e}om.\ Diff.} 27 (1986), 
3-17. Corrections in \cite{KR2}.
\bibitem{CCBI} A.\ Kock, Combinatorics of curvature and the Bianchi 
identity, {\em Theory and Appl.\ of Categories} 2 (1996), 69-89.

\bibitem{PBGC} A.\ Kock, Principal bundles, groupoids, and 
connections, in {\em Geometry and Topology of Manifolds} 
(ed.\ J.\ Kubarski, J.\ Pradines, T.\ Rybicki and R.\ Wolak), Banach Center 
Publications Vol.\ 76 (2007), 185-200.





\bibitem{END}A.\ Kock: Envelopes - notion and definiteness, 
{\em Beitr\"{a}ge  zur Algebra und Geometrie} 48 (2007), 345-350.
 \bibitem{SGM}A.\ Kock, {\em Synthetic Geometry of Manifolds}, 
Cambridge Tracts in Mathema\-tics 180, Cambridge University Press 2010.
\bibitem{MSSDG} A.\ Kock, Metric spaces and SDG, {\em Theory and Appl.\ of Categories} 
32 (2017), 803-822.
\bibitem{KM}A.\ Kock and C.J.\ Mikkelsen, Topos theoretic 
factorization of non-standard extensions, in {\em Victoria Symposium 
on Nonstandard Analysis 1972} (ed.\ A.\ Hurd and P.\ Loeb),   
Springer Lecture Notes in Math.\ 369 (1974), 122-143. 
\bibitem{KR2}A.\ Kock and G.E.\ Reyes, Corrigendum and addenda to 
``Convenient vector spaces embed'', {\em Cahiers de Top.\ et 
G\'{e}om.\ Diff.} 28 (1987), 69-89.
\bibitem{KR3}A.\ Kock and G.E.\ Reyes,  Aspects of fractional 
exponent functors, {\em Theory and Appl.\ of Categories} 5 (1999), 251-265.
\bibitem{KumS}A.\ Kumpera and D.\ Spencer, {\em Lie equations}, Annals of 
Mathematics Studies 73, Princeton 1972.
\bibitem{Lav}R.\ Lavendhomme, {\em Basic Concepts Of Synthetic
Differential Geometry}, Kluwer Academic Publishers 1996.
\bibitem{L98}F.W.\ Lawvere, Outline of Synthetic Differential 
Geometry,  Notes Buffalo 1998,
http://www.acsu.buffalo.edu/wlawvere/SDG\_  Outline.pdf
\bibitem{LC} F.W.\ Lawvere, Axiomatic cohesion, {\em Theory  and 
Appl. of Categories} 
19 (2007), 41-47.
 \bibitem{L:D}F.W.\ Lawvere, Euler's Continuum Functorially 
Vindicated, vol.\ 75 of {\em The Western Ontario Series in Philosophy of Science}, 2011.
\bibitem{LaMe} F.W.\ Lawvere and M.\ Menni, Internal choice holds in 
the discrete part of any cohesive topos satisfying stable connected 
codiscreteness, {\em Theory and Appl.\ of Categories} 30 (2015), 
909-932. 
\bibitem{Lie} S.\ Lie, Allgemeine Theorie der partiellen 
Differentialgleichungen erster Ordnung, {\em Math.\ Ann.\ } 9 (1876), 245-296.
\bibitem{BT} S.\ Lie, {\em Geometrie der Ber\"{u}hrungstransformationen}, Leipzig 
1896, reprint Chelsea Publ.\ Comp.\ 1977.
\bibitem{May}P.\  May, {\em Simplicial Objects in Algebraic 
Topology}, van Nostrand Math.\ Studies 11, 1967.
\bibitem{Mcl}C.\ McLarty, Local, and some global, results in 
synthetic differential geometry, in {\em Category Theoretic Methods in 
Geometry} (ed.\ A.\ Kock), Aarhus Mat.\ Inst.\ Various Publ.\ Series 35 
(1983), 226-256.
\bibitem{Menni} M.\ Menni, Continuous cohesion over sets, 
{\em Theory and Appl.\ of Categories} 29 (2014), 542-568.
\bibitem{MR}I.\ Moerdijk  and
G.E.\ Reyes, {\em Models for Smooth Infinitesimal Analysis}, Springer 1991.
\bibitem{Mumford} D.\ Mumford, {\em The Red Book of Varieties and 
Schemes},  $\leq$ 1968, reprinted 1988 as Springer Lecture Notes in Math.\ 
1358. 
\bibitem{Nelson}  E.\ Nelson, {\em Internal set theory: A new approach to nonstandard analysis}, Bull. Amer.
Math. Soc. 83 (1977), 1165-1198.
\bibitem{N}H.\ Nishimura, Higher-Order Preconnections in 
Synthetic Differential Geometry of Jet Bundles, {\em Beitr\"{a}ge zur Algebra 
und Geometrie} 45 (2004), 677-696.
\bibitem{Weil}A.\ Weil, Th\'{e}orie des points proches sur les vari\'{e}t\'{e}s 
diff\'{e}rentiables, in {\em Colloq.\  Top.\  et G\'{e}om.\ Diff.\ }, Strassbourg 
1953.
\end{thebibliography}
\end{document}